\documentclass[a4paper]{article}

\usepackage{amsmath,amssymb}
\usepackage{amsthm}
\usepackage{a4wide}
\usepackage{graphicx,color,tikz}
\newcommand{\sign}{\text{sign}}

\newtheorem{proposition1}{Proposition}
\newtheorem{theorem}{Theorem}
\newtheorem{remark}{Remark}

\begin{document}

\begin{center}

\begin{Large}
\textbf{Optimized Schwarz methods for the time-dependent Stokes-Darcy coupling}
\end{Large}

\medskip

Marco Discacciati$^{1}$ and Tommaso Vanzan$^{2}$

\medskip

${}^1$ Department of Mathematical Sciences, Loughborough University, Loughborough LE11 3TU, United Kingdom, m.discacciati@lboro.ac.uk.

${}^2$ D\'epartement de Math\'ematiques, \'Ecole Polytechnique F\'ed\'erale de Lausanne, CH-1015 Lausanne, Switzerland, tommaso.vanzan@epfl.ch.

\end{center}

\bigskip

\centerline{\textbf{Abstract}}


This paper derives optimal coefficients for optimized Schwarz iterations for the time-dependent Stokes-Darcy problem using an innovative strategy to solve a nonstandard min-max problem. The coefficients take into account both physical and discretization parameters that characterize the coupled problem, and they guarantee the robustness of the associated domain decomposition method. Numerical results validate the proposed approach in several test cases with physically relevant parameters.

\bigskip

\section{Introduction}

The Stokes-Darcy problem has been extensively studied during the last two decades due to its relevance to model filtration phenomena in industrial and natural applications. The steady problem introduced in the seminal works \cite{Discacciati:2002:MNM,Layton:2003:CFF} has been extended to consider the time dependent case in, e.g., \cite{Discacciati:2004:IMS,Cao:2010:FEA,Cesmelioglu:2013:TDC,Cao:2014:PNI,Moraiti:2012:JMAA,Rybak:2014:JCP}.

In both settings, the space discretization of the Stokes-Darcy model leads to a large linear system (possibly at each time step) that has to be effectively preconditioned. A possible strategy to achieve this is to adopt a monolithic approach considering the whole coupled linear system at once. (see, e.g., the recently proposed robust monolithic preconditioner derived using an operator preconditioning framework in \cite{Mardal:2011:PDS}). However, the multi-physics nature of the Stokes-Darcy problem makes it suitable for  
%
%
decoupled strategies based on domain decomposition that set up an iterative process where the Stokes and Darcy problems are solved separately at each iteration until convergence. 
Nevertheless, the computational efficiency of a domain decomposition approach depends on the number of subdomain iterations that are needed to reach convergence, especially in the time-dependent case. It is now well understood that they way in which the subdomain are coupled at each iteration significantly affects the convergence rate.

Earlier works focused on Dirichlet-Neumann algorithms \cite{Discacciati:2004:CAS,Discacciati:2004:PHD,Discacciati:2009:RMC}, which however exhibit slow convergence for small values of the viscosity of the fluid and the permeability of the porous medium. 
More recent efforts have focused on Robin-Robin transmission conditions \cite{Discacciati:2007:RRD,Chen:2011:PRR,Cao:2011:RRD,Caiazzo:2014:CIS} which generally show better properties in terms of convergence and robustness with respect to the physical parameters, provided that the Robin parameters are properly selected. Their optimization is usually carried out in a simplified geometrical setting using Fourier analysis. Robin-Robin domain decomposition methods that use optimized Robin parameters belongs to the family of optimized Schwarz methods \cite{Gander:2006:OSM}, which have been proven to be very effective, beyond the limitations set by the Fourier analysis, for several different equations \cite{LGG:SISC:2009,Gander_XU:2016:OSM,gander2019heterogeneous,Gander:2019:SIREV} and geometric configurations \cite{Chen:2021:OSM,Gigante:2021:OSM,Gigante:2020:OSM}.
Optimized Schwarz methods have been first studied for the stationary Stokes-Darcy system in \cite{Discacciati:2018:IMAJNA} and \cite{Gander:2020,gander:MOSM}.

Few works have considered the time-dependent Stokes-Darcy coupling with Robin interface conditions. References \cite{Feng:2012:AMC,Cao:2014:PNI} present a loosely coupled implicit time-marching scheme that at each time-step solves the Stokes and Darcy problems separately, without iterating. The two subdomains are coupled through a Robin boundary condition that depends on the solution computed in the other subdomain at the previous time step. Despite being non-iterative and preserving optimal accuracy, this approach has the downside that the iterates do not fulfill the physical coupling conditions at each time step. Time parallel strategies based on waveform relaxation algorithms have instead being analyzed in \cite{Hoang:2022:NCT}.

In this work, we propose a novel approach where we perform a classical semi-discretization in time using a general implicit $\theta$-method \cite{wanner2002geometric}, and we optimize the convergence rate of a non-overlapping Schwarz methods to solve the local subproblems at each time step. For the optimization of the Schwarz method, we rely on Fourier analysis to derive the convergence factor of the iterative domain decomposition method in the frequency domain, and we study a min-max problem to characterize the optimized Robin parameters. Since the min-max problem differs from others in the existing literature, we develop a new theoretical argument that can be applied to min-max sharing the same structure.
Numerical results performed with physically relevant parameters show that the optimized Schwarz method is robust and requires very few iterations to convergence at each time-step, so that the method we propose is a valid alternative to non-iterative time-marching schemes.

The rest of the manuscript is organized as follows. Section \ref{sec:Setting} formulates the problem, describes the semi-discretization in time and introduces the optimized Schwarz method. Section \ref{sec:Robin-Robin} performs a Fourier analysis and derives the optimized parameters by solving the min-max problem. Section \ref{sect:NumericalResults} presents extensive numerical tests to assess the performance of the proposed algorithm.

\section{Setting and formulation of the optimized Schwarz method}\label{sec:Setting}

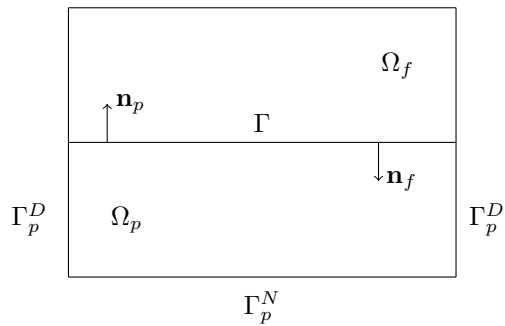
\begin{figure}[bht]
\begin{center}
\begin{tikzpicture}[scale=0.51]
\draw (1,0.5)--(11,0.5);
\draw (1,0.5)--(1,7.5);
\draw (1,7.5)--(11,7.5);
\draw (11,0.5)--(11,7.5);
\draw (1,4)--(11,4);
\draw [->] (2,4) -- (2,5); 
\node at (2.6,5) {$\mathbf{n}_p$};
\draw [->] (9,4) -- (9,3); 
\node at (9.6,3) {$\mathbf{n}_f$};
\node at (6,4.5) {$\Gamma$};
\node at (9.5,6) {$\Omega_f$};
\node at (2.5,2) {$\Omega_p$};
\node at (6,-0.3) {$\Gamma_p^N$};
\node at (11.8,2) {$\Gamma_p^D$};
\node at (0,2) {$\Gamma_p^D$};
\end{tikzpicture}
\caption{Schematic representation of a 2D section of the computational domain.}
\label{fig:domain}
\end{center}
\end{figure}

Let $\Omega \subset \mathbb{R}^\mathrm{D}$ ($\mathrm{D}=2,3$) be an open bounded domain formed by a fluid region $\Omega_f$ and a porous-medium region $\Omega_p$ with Lipschitz continuous boundaries $\partial \Omega_f$ and $\partial \Omega_p$. The two regions are non-overlapping and separated by an interface $\Gamma$, i.e., $\overline{\Omega} = \overline{\Omega}_f \cup \overline{\Omega}_p$, $\Omega_f \cap \Omega_p = \emptyset$ and $\overline{\Omega}_f \cap \overline{\Omega}_p = \Gamma$, as illustrated in Fig.~\ref{fig:domain}. Let $\mathbf{n}_p$ and $\mathbf{n}_f$ be the unit normal vectors pointing outwards of $\partial \Omega_p$ and $\partial \Omega_f$, respectively, with $\mathbf{n}_f=-\mathbf{n}_p$ on $\Gamma$. We assume that $\mathbf{n}_f$ and $\mathbf{n}_p$ are regular enough, and we indicate $\mathbf{n} = \mathbf{n}_f$ for simplicity.


In $\Omega_f$ and for time $t \in (0,T]$, we consider an incompressible fluid with constant viscosity and density described by the dimensionless time-dependent Stokes equations:
\begin{eqnarray}
\partial_t \mathbf{u}_f - \boldsymbol{\nabla} \cdot (2\mu_f \nabla^{\mathrm{s}} \mathbf{u}_f - p_f \boldsymbol{I}) = \mathbf{f}_f && \mbox{in } \Omega_f\,, \label{eq:stokes1} \\
\nabla \cdot \mathbf{u}_f = 0 && \mbox{in } \Omega_f\,\nonumber,
\end{eqnarray}
where $\mu_f = Re^{-1}$, $Re$ being the Reynolds number, $\mathbf{u}_f$ and $p_f$ are the fluid velocity and pressure, $\boldsymbol{I}$ and $\nabla^{\mathrm{s}} \mathbf{u}_f = \frac{1}{2} (\nabla \mathbf{u}_f + (\nabla \mathbf{u}_f)^T)$ are the identity and the strain rate tensor, and $\mathbf{f}_f$ is an external force.
In the porous medium domain $\Omega_p$ and for $t \in (0,T]$, we consider the dimensionless time-dependent Darcy model (see, e.g., \cite{Cao:2010:FEA,Moraiti:2012:JMAA,Rybak:2014:JCP}):
\begin{equation}\label{eq:darcy}
S_p \, \partial_t \,p_p - \nabla \cdot (\boldsymbol{\eta}_p \nabla p_p) = f_p \qquad \mbox{in } \Omega_p\,,
\end{equation}
where $p_p$ is the fluid pressure in the porous medium, $\boldsymbol{\eta}_p$ is the permeability tensor, $S_p$ is the specific mass storativity coefficient, and $f_p$ is an external force. 

The two models are coupled through the classical Beaver-Joseph-Saffman conditions which describe the filtration across the interface $\Gamma$ \cite{Beavers:1967:BCN,Jager:1996:OBC,Rybak:2014:JCP,Saffman:1971:OBC}:
\begin{eqnarray}
\mathbf{u}_f \cdot \mathbf{n} = - ( \boldsymbol{\eta}_p \nabla p_p ) \cdot \mathbf{n} && \mbox{at } \Gamma, \label{eq:interf1} \\
- \mathbf{n} \cdot ( 2\mu_f \nabla^{\mathrm{s}} \mathbf{u}_f - p_f \boldsymbol{I} ) \cdot \mathbf{n} = p_p && \mbox{at } \Gamma, \label{eq:interf2} \\
- ( ( 2\mu_f \nabla^{\mathrm{s}} \mathbf{u}_f - p_f \boldsymbol{I} ) \cdot \mathbf{n})_\tau = \xi_f (\mathbf{u}_f)_\tau  && \mbox{at } \Gamma, \label{eq:interf3}
\end{eqnarray}
where $\xi_f = \alpha_{BJ} \mu_f/\sqrt{{\boldsymbol{\tau}\cdot \boldsymbol{\eta}_p \cdot \boldsymbol{\tau}}}$ and $\alpha_{BJ}$ is a dimensionless constant which depends on the geometric structure of the porous medium. The notation $(\mathbf{v})_{\tau}$ indicates the tangential component of any vector $\mathbf{v}$ at $\Gamma$, i.e., $(\mathbf{v})_\tau = \mathbf{v} - ( \mathbf{v} \cdot \mathbf{n} ) \, \mathbf{n}$ at $\Gamma$.
As boundary conditions, we impose $p_p = 0$ on $\Gamma_p^D$, $\mathbf{u}_p \cdot \mathbf{n}_p = 0$ on $\Gamma_p^N$, and $\mathbf{u}_f=\mathbf{0}$ on $\partial\Omega_f \setminus \Gamma$ as in \cite{Discacciati:2018:IMAJNA} (see Fig. \ref{fig:domain} for the notation), and as initial conditions we set
\begin{equation}\label{eq:initialCond}
\mathbf{u}_f = \mathbf{u}_f^0 \quad \mbox{in } \Omega_f
\quad \mbox{and} \quad
p_p=p_p^0 \quad \mbox{in } \Omega_p \quad \mbox{at } t=0,
\end{equation}
with $\mathbf{u}_f^0$ and $p_p^0$ given functions.

We define the spaces
\begin{equation*}
\begin{array}{c}
\mathbf{V}_f := \{ \mathbf{v} \in \mathbf{H}^1(\Omega_f) \, : \, \mathbf{v} = \mathbf{0} \quad \mbox{on } \partial\Omega_f \setminus \Gamma \, \}, \quad
Q_f := L^2(\Omega_f) \, ,\\[5pt]
Q_p := \{ q \in H^1(\Omega_p) \, : \, q=0 \quad \mbox{on } \Gamma_p^D \, \} \, ,
\end{array}
\end{equation*}
and let $(\cdot,\cdot)_{\mathcal{D}}$ be the $L^2$ scalar product in a domain $\mathcal{D}$ for scalar, vector, or tensor functions, while $\langle \cdot, \cdot \rangle_\Gamma$ is the scalar product in $H^{1/2}(\Gamma)$. Then, we introduce the bilinear forms
\begin{eqnarray*}
a_f : \mathbf{V}_f \times \mathbf{V}_f \to \mathbb{R},
 && a_f (\mathbf{v},\mathbf{w}) = ( 2\mu_f \nabla^{\mathrm{s}}\mathbf{v} ,\nabla^{\mathrm{s}}\mathbf{w} )_{\Omega_f} + \xi_f \langle (\mathbf{v})_\tau, (\mathbf{w})_\tau \rangle_\Gamma , \\
b_f: \mathbf{V}_f \times Q_f\rightarrow \mathbb{R},
 && b_f(\mathbf{v},q)= -(\nabla \cdot \mathbf{v} ,q)_{\Omega_f}, \\
a_p : \, Q_p \times Q_p \to \mathbb{R},
 && a_p (p,q) = \left( \boldsymbol{\eta}_p \, \nabla p , \nabla q \right)_{\Omega_p} .
\end{eqnarray*}

The weak formulation of problem \eqref{eq:stokes1}--\eqref{eq:interf3} can be written as: for all $t \in (0,T]$, find $\mathbf{u}_f(t) \in \mathbf{V}_f$, $p_f(t) \in Q_f$, and $p_p(t) \in Q_p$ such that
\begin{equation}\label{eq:weakForm}
\begin{array}{rcll}
( \partial_t \mathbf{u}_f , \mathbf{v}_f )_{\Omega_f} + a_f (\mathbf{u}_f,\mathbf{v}_f) + b_f (\mathbf{v}_f,p_f) + \langle p_p, \mathbf{v}_f \cdot \mathbf{n} \rangle_\Gamma &=& ( \mathbf{f}_f , \mathbf{v} )_{\Omega_f}
 & \quad \forall \mathbf{v}_f \in \mathbf{V}_f, \\[5pt]
b_f(\mathbf{u}_f, q_f) &=& 0
 & \quad \forall q_f \in Q_f, \\[5pt]
( S_p\, \partial_t \,p_p , q_p )_{\Omega_p} + a_p (p_p, q_p) - \langle \mathbf{u}_f \cdot \mathbf{n}, q_p \rangle_\Gamma &=& ( f_p , q )_{\Omega_p}
 & \quad \forall q_p \in Q_p,
\end{array}
\end{equation}
with initial conditions \eqref{eq:initialCond}.
For all $t\in (0,T]$, the weak formulation \eqref{eq:weakForm} admits a unique solution $(\mathbf{u}_f(t),p_f(t),p_p(t))\in \mathbf{V}_f\times Q_f\times Q_p$ which depends continuously on the data (see, e.g., \cite{Discacciati:2004:PHD}). 
We remark that the scalar product $\langle \cdot , \cdot \rangle_\Gamma$ can be used both in \eqref{eq:weakForm} and in the definition of the bilinear form $a_f$ due to the regularity assumption on $\mathbf{n}$ at $\Gamma$.

\subsection{Optimized Schwarz method}

A linear combination of the interface conditions \eqref{eq:interf1} and \eqref{eq:interf2} with coefficients $(-\alpha_f,1)$ and $(\alpha_p,1)$, $\alpha_f,\alpha_p$ being positive real parameters, leads to the Robin transmission conditions
\begin{equation}\label{eq:interfrobin1}
- \mathbf{n} \cdot ( 2\mu_f \nabla^{\mathrm{s}} \mathbf{u}_f - p_f \boldsymbol{I} ) \cdot \mathbf{n} -\alpha_f \mathbf{u}_f \cdot \mathbf{n} =  p_p + \alpha_f \, \boldsymbol{\eta}_p \nabla p_p\cdot \mathbf{n} \quad \mbox{at } \Gamma\,,
\end{equation}
and
\begin{equation}\label{eq:interfrobin2}
 p_p - \alpha_p \, \boldsymbol{\eta}_p \nabla p_p\cdot \mathbf{n} = - \mathbf{n} \cdot ( 2\mu_f \nabla^{\mathrm{s}} \mathbf{u}_f - p_f \boldsymbol{I} ) \cdot \mathbf{n} + \alpha_p \mathbf{u}_f \cdot \mathbf{n} \quad \mbox{at } \Gamma \, .
\end{equation}

Upon setting
\begin{equation*}
\lambda_{f}(t)=p_p+\alpha_f \,\boldsymbol{\eta}_p \nabla p_p\cdot \mathbf{n}
\quad \mbox{and} \quad
\lambda_p(t)= - \mathbf{n} \cdot ( 2\mu_f \nabla^{\mathrm{s}} \mathbf{u}_f - p_f \boldsymbol{I} )\cdot \mathbf{n} + \alpha_p \, \mathbf{u}_f\cdot \mathbf{n} \, ,
\end{equation*}
the Robin interface conditions \eqref{eq:interfrobin1} and \eqref{eq:interfrobin2} can be used to equivalently reformulate problem \eqref{eq:weakForm} as: for all $t\in (0,T]$, find $\mathbf{u}_f(t) \in \mathbf{V}_f$, $p_f(t) \in Q_f$, and $p_p(t) \in Q_p$ such that
\begin{equation}\label{eq:weakForm_Robin}
\hspace*{-3mm}
\begin{array}{rcl}
( \partial_t \mathbf{u}_f , \mathbf{v}_f )_{\Omega_f} + a_f (\mathbf{u}_f,\mathbf{v}_f) + b_f (\mathbf{v}_f,p_f) + \alpha_f \langle \mathbf{u}_f\cdot \mathbf{n}, \mathbf{v}_f \cdot \mathbf{n} \rangle_\Gamma &=& ( \mathbf{f}_f , \mathbf{v} )_{\Omega_f} -\langle \lambda_f, \mathbf{v}_f\cdot \mathbf{n}\rangle_\Gamma
 \\[5pt]
b_f(\mathbf{u}_f, q_f) &=& 0 \\[5pt]
( S_p\, \partial_t \,p_p , q_p )_{\Omega_p} + a_p (p_p, q_p) + \frac{1}{\alpha_p}\langle p_p, q_p \rangle_\Gamma &=& ( f_p , q )_{\Omega_p} +\frac{1}{\alpha_p}\langle \lambda_{p}, q_p\rangle_{\Gamma},
\end{array}
\end{equation}
for all $\mathbf{v}_f \in \mathbf{V}_f$, $q_f \in Q_f$, and $q_p \in Q_p$.

For the time discretization of problem \eqref{eq:weakForm_Robin}, we split the interval $[0,T]$ into equal subintervals $[t^{n-1},t^n]$ such that $[0,T]=\cup_{n=1}^N[t^{n-1},t^n]$, with $t^0=0$, $t^N=T$ and $\Delta t=t^{n}-t^{n-1}$, and we consider the $\theta$-method (see, e.g., \cite{wanner2002geometric,Turek:1996:IJNMF,John:2006:CMAME}) with $\theta \in (0,1]$. We are particularly interested in the cases $\theta=1$ and $\theta=\frac{1}{2}$ which correspond to the implicit Euler method and to the Crank-Nicolson method.

\medskip

For simplicity of notation, let us define the bilinear forms
\begin{eqnarray*}
\widetilde{a}_f : \mathbf{V}_f \times \mathbf{V}_f \to \mathbb{R},
 && \widetilde{a}_f (\mathbf{v},\mathbf{w}) = (\mathbf{v},\mathbf{w})_{\Omega_f} + \theta \, \Delta t \left( a_f(\mathbf{v},\mathbf{w}) +\alpha_f \langle \mathbf{v}\cdot \mathbf{n},\mathbf{w}\cdot \mathbf{n}\rangle_{\Gamma} \right) \,,\\
\widetilde{a}_p :\, Q_p \times Q_p \to \mathbb{R},
 && \widetilde{a}_p (p,q) = (S_p\, p,q)_{\Omega_p} + \theta \, \Delta t \left( a_p(p ,q) +\frac{1}{\alpha_p}\langle p,q\rangle_\Gamma\right) ,
\end{eqnarray*}
and the functionals $F_f^n:\mathbf{V}_f\to \mathbb{R}$ and ${F}_p^n:Q_p\to \mathbb{R}$
\begin{equation*}
\begin{array}{rcl}
{F}_f^{n}(\mathbf{v}) &=& \Delta t \left( \, \theta \, (\mathbf{f}_f^n,\mathbf{v})_{\Omega_f} + (1-\theta) (\mathbf{f}_f^{n-1},\mathbf{v})_{\Omega_f} \right) + (\mathbf{u}^{n-1}_f,\mathbf{v})_{\Omega_f} \\
&& -\,(1-\theta)\,\Delta t\left(a_f(\mathbf{u}^{n-1}_f,\mathbf{v})+\alpha\langle \mathbf{u}_f^{n-1}\cdot \mathbf{n},\mathbf{v}\cdot \mathbf{n}\rangle +\langle \lambda_f^{n-1},\mathbf{v}_f\cdot \mathbf{n}\rangle \right),\\[5pt]
{F}_p^{n}(q) &=& \Delta t \left( \, \theta \, (f_p^n,q)_{\Omega_p} + (1-\theta) (f_p^{n-1},q)_{\Omega_p} \right)  + ( S_p\,  p^{n-1}_p , q )_{\Omega_p} \\[3pt]
 && - (1-\theta)\,\Delta t \, \left(a_p(p_p^{n-1},q)+\frac{1}{\alpha_p}\langle p_p^{n-1},q\rangle_\Gamma -\frac{1}{\alpha_p}\langle \lambda_p^{n-1},q_p\rangle\right).
\end{array}
\end{equation*}

Then, starting from the initial condition \eqref{eq:initialCond}, for $n=1,\dots,N$, we compute $(\mathbf{u}_f^n,p_f^n,p_p^n)\in \mathbf{V}_f \times Q_f \times Q_p$ solution of
\begin{equation}\label{eq:weakForm_Robinsemidisc_v2}
\begin{array}{rcll}
\widetilde{a}(\mathbf{u}^n_f,\mathbf{v}_f^n) + \Delta t \, b_f (\mathbf{v}_f,p^{n}_f) &=& F_f^{n}(\mathbf{v}_f)-\theta\,\Delta t\,\langle \lambda^n_f, \mathbf{v}_f\cdot \mathbf{n}\rangle_\Gamma
 & \quad \forall \mathbf{v}_f \in V_f, \\[5pt]
b_f(\mathbf{u}^{n}_f, q_f) &=& 0
 & \quad \forall q_f \in Q_f, \\[5pt]
\widetilde{a}_p(p_p^n,q_p) &=& {F}^{n}_p(q_p) + \frac{\theta \, \Delta t}{\alpha_p}\langle \lambda^n_p, q_{p}\rangle_{\Gamma}
 & \quad \forall q_p \in Q_p,
\end{array}
\end{equation}
where
\begin{equation*}
\lambda^n_{f}=p^n_p+\alpha_f \,\boldsymbol{\eta}_p \nabla p^n_p\cdot \mathbf{n},
\qquad
\lambda^n_p= - \mathbf{n} \cdot ( 2\mu_f \nabla^{\mathrm{s}} \mathbf{u}^n_f - p^n_f \boldsymbol{I} )\cdot \mathbf{n} +\alpha_p \,\mathbf{u}^n_f\cdot \mathbf{n} \, .
\end{equation*}

At each time level $n\geq 1$, system \eqref{eq:weakForm_Robinsemidisc_v2} is still fully coupled. To decouple the Stokes and the Darcy equations, we consider the Robin-Robin iterative scheme which starts from the initial guess $\lambda_f^{n,0}$, 
and, for $m\geq 1$ until convergence, it computes 
\begin{equation}\label{eq:weakForm_Robinsemidisc_iter_v2}
\begin{array}{rcll}
\widetilde{a}(\mathbf{u}^{n,m}_f,\mathbf{v}_f^n) + \Delta t \, b_f (\mathbf{v}_f,p^{n,m}_f) &=& F_f^{n}(\mathbf{v}_f)-\theta\,\Delta t\,\langle \lambda^{n,m-1}_f, \mathbf{v}_f\cdot \mathbf{n}\rangle_\Gamma
 & \quad \forall \mathbf{v}_f \in V_f, \\[5pt]
b_f(\mathbf{u}^{n,m}_f, q_f) &=& 0
 & \quad \forall q_f \in Q_f, \\[5pt]
\widetilde{a}_p(p_p^{n,m},q_p) &=& {F}^{n}_p(q_p) + \frac{\theta \, \Delta t}{\alpha_p}\langle \lambda^{n,m}_p, q_{p}\rangle_{\Gamma}
 & \quad \forall q_p \in Q_p,
\end{array}
\end{equation}
where
\begin{equation*}
\lambda^{n,m}_{f}=p^{n,m}_p+\alpha_f \,\boldsymbol{\eta}_p \nabla p^{n,m}_p\cdot \mathbf{n},
\quad
\lambda^{n,m}_p= - \mathbf{n} \cdot ( 2\mu_f \nabla^{\mathrm{s}} \mathbf{u}^{n,m}_f - p^{n,m}_f \boldsymbol{I} )\cdot \mathbf{n} +\alpha_p \,\mathbf{u}^{n,m}_f\cdot \mathbf{n} \, .
\end{equation*}

Standard calculations (see, e.g., \cite{Discacciati:2018:IMAJNA,Gander:2020}) show that, for an initial guess $\lambda_p^{n,0}$, for $m \geq 1$,
\begin{equation}\label{eq:relationlambda_v2}
\begin{aligned}
\lambda^{n,m}_f &= \left(1+\frac{\alpha_f}{\alpha_p}\right)p_p^{n,m} -\frac{\alpha_f}{\alpha_p}\lambda_{p}^{n,m-1},\\
\lambda^{n,m}_p &= (\alpha_f+\alpha_p) \, \mathbf{u}_f^{n,m} \cdot \mathbf{n} +\lambda_f^{n,m}.
\end{aligned}
\end{equation}
Notice that a parallel version of the Robin-Robin scheme can be obtained by taking $\lambda_p^{n,m-1}$ in \eqref{eq:weakForm_Robinsemidisc_iter_v2}$_3$ and $\lambda_f^{n,m-1}$ in \eqref{eq:relationlambda_v2}$_2$.

\medskip

We conclude this section by showing that the iteration scheme \eqref{eq:weakForm_Robinsemidisc_iter_v2} can be reformulated only in terms of the interface functions $\lambda_{f}^{n,m}$ and $\lambda_p^{n,m}$ at $\Gamma$. To this aim, in the domain $\Omega_f$, we define the continuous trace operator $\tau_f:\mathbf{V}_f\times Q_f\rightarrow H^{\frac{1}{2}}(\Gamma)$ as $\tau_f((\mathbf{v},p))=(\mathbf{v}\cdot \mathbf{n})_{|\Gamma}$,  and the continuous extension operator $\mathcal{E}_f: H^{-\frac{1}{2}}(\Gamma)\times \mathbf{V}_f^{-1}\rightarrow \mathbf{V}_f \times Q_f$ as $\mathcal{E}_f(\lambda,\mathcal{F}_f)=(\mathbf{v},p)$, where $(\mathbf{v},p)$ is the solution of 
\begin{equation*}
\begin{array}{rcll}
\widetilde{a}_f(\mathbf{v},\mathbf{w}) + \Delta t \, b_f (\mathbf{v},p) &=&\mathcal{F}_f(\mathbf{w}) - \theta \, \Delta t \, \langle \lambda, \mathbf{w}\cdot \mathbf{n}\rangle_\Gamma
 & \quad \forall \mathbf{w} \in \mathbf{V}_f, \\[5pt]
b_f(\mathbf{v}, q) &=& 0
 & \quad \forall q \in Q_f. \\[5pt]
\end{array}
\end{equation*}
Similarly in the domain $\Omega_p$, we define the continuous trace operator $\tau_p:Q_p\rightarrow H^{\frac{1}{2}}(\Gamma)$ as $\tau_p(q)=q_{|\Gamma}$, and the continuous extension operator
$\mathcal{E}_p: H^{-\frac{1}{2}}(\Gamma)\times Q_p^{-1}\rightarrow Q_p$ as $\mathcal{E}_p(\lambda,\mathcal{F}_p)=p$, where $p$ is the solution of 
\begin{equation*}
\widetilde{a}_p(p,q) = \mathcal{F}_p(q)+\frac{\theta \,\Delta t}{\alpha_p}\langle \lambda,q\rangle
 \qquad \forall q \in Q_p \,.
\end{equation*}

Then, relations \eqref{eq:relationlambda_v2} can be equivalently reformulated as
\begin{equation}\label{eq:relationlambda2_v2}
\begin{aligned}
\lambda^{n,m}_f &= \left( 1+\frac{\alpha_f}{\alpha_p} \right) \, \tau_p(\mathcal{E}_p(\lambda^{n,m-1}_p, {F}_p^{n}))-\frac{\alpha_f}{\alpha_p} \lambda_p^{n,m-1},\\[5pt]
\lambda^{n,m}_p &= (\alpha_f+\alpha_p) \, \tau_f(\mathcal{E}_{f}(\lambda^{n,m}_f,{F}_f^{n}))+\lambda_f^{n,m}.
\end{aligned}
\end{equation}

Now define the linear continuous operators
\begin{equation*}
\begin{array}{ll}
\mathcal{G}_f: H^{-\frac{1}{2}}(\Gamma) \rightarrow H^{\frac{1}{2}}(\Gamma), \qquad \mathcal{G}_f(\lambda):=\tau_f\left(\mathcal{E}_f(\lambda,0)\right)\,,\\[4pt]
\mathcal{G}_p: H^{-\frac{1}{2}}(\Gamma) \rightarrow H^{\frac{1}{2}}(\Gamma), \qquad
\mathcal{G}_p(\lambda):=\tau_p\left(\mathcal{E}_p(\lambda,0)\right)\,.
\end{array}
\end{equation*}
Thanks to the linearity of the extension operators $\mathcal{E}_f$ and $\mathcal{E}_p$, and setting
\begin{equation*}
\chi^{n}_f= \left( 1+\frac{\alpha_f}{\alpha_p} \right) \tau_{p} (\mathcal{E}_{p}(0,{F}_{p}^{n}))\,,
\qquad
\chi^{n}_p= (\alpha_f + \alpha_p) \, \tau_{f}(\mathcal{E}_{f}(0,{F}_{f}^{n}))\,,
\end{equation*}
equation \eqref{eq:relationlambda2_v2} can be equivalently written as 
\begin{equation}\label{eq:stationary_iteration}
\begin{pmatrix}
\lambda_f^{n,m}\\\lambda_p^{n,m}
\end{pmatrix}=\begin{pmatrix}
0 & -\frac{\alpha_f}{\alpha_p}I +(1+\frac{\alpha_f}{\alpha_p})\, \mathcal{G}_{p}(\cdot)\\
(\alpha_f+\alpha_p) \, \mathcal{G}_{f}(\cdot) + I & 0
\end{pmatrix}\begin{pmatrix}
\lambda_f^{n,m-1}\\\lambda_p^{n,m\phantom{-1}}
\end{pmatrix}+\begin{pmatrix}
\chi^{n}_f\\ \chi^{n}_p
\end{pmatrix}.
\end{equation}

This corresponds to a Gauss-Seidel step for the linear interface system
\begin{equation}\label{eq:interfaceSystem_v2}
\begin{pmatrix}
I & \frac{\alpha_f}{\alpha_p}I - (1+\frac{\alpha_f}{\alpha_p})\, \mathcal{G}_{p}(\cdot)\\
-(\alpha_f+\alpha_p) \, \mathcal{G}_{f}(\cdot) - I & I
\end{pmatrix}\begin{pmatrix}
\lambda_f^{n}\\\lambda_p^{n}
\end{pmatrix}
=
\begin{pmatrix}
\chi^{n}_f\\ \chi^{n}_p
\end{pmatrix}
\end{equation}
which can be more effectively solved using a suitable Krylov method (e.g., GMRES \cite{Saad:1986:GMR}) as we will do in Sect.~\ref{sect:NumericalResults}.

\section{Analysis of the optimized Schwarz method}\label{sec:Robin-Robin}

In this section, we analyse the optimized Schwarz method \eqref{eq:weakForm_Robinsemidisc_v2} with the aim of characterizing optimal values of the parameters $\alpha_f$ and $\alpha_p$ to accelerate the convergence of the method. The expression of the convergence rate is obtained in Sect.~\ref{sect:ConvergenceRate}, while the optimization of the parameters is carried out in Sect.~\ref{sect:Optimization} under the simplifying hypothesis that the Beavers-Joseph-Saffman interface condition \eqref{eq:interf3} is replaced by the zero tangential velocity condition $(\mathbf{u}_f)_\tau = 0$ at $\Gamma$. 

\subsection{Characterization of the convergence rate}\label{sect:ConvergenceRate}

Since the problem is linear, for the analysis of the optimized Schwarz method \eqref{eq:weakForm_Robinsemidisc_v2} we directly consider the error equation obtained by setting $\mathbf{f}_f=\mathbf{0}$ and $f_p=0$. More precisely, at time step $n \geq 1$, for $m\geq 1$ until convergence, we consider
\begin{enumerate}
\item
the Stokes problem: find $\mathbf{u}_f^{n,m}$, $p_f^{n,m}$ such that
\begin{equation}\label{eq:Stokes_problem_v2}
\hspace*{-3mm}
\begin{array}{l l l}
\displaystyle \frac{\mathbf{u}_f^{n,m}}{\Delta t} - \theta\,\boldsymbol{\nabla} \cdot ( 2\mu_f\, \nabla^{\mathrm{s}} \mathbf{u}_f^{n,m}) + \nabla p_f^{n,m} = \mathbf{0} && \mbox{in } \Omega_f\,, \\[5pt]
\nabla \cdot \mathbf{u}_f^{n,m} = 0 && \mbox{in } \Omega_f\,, \\[3pt]
- ( ( 2\mu_f \nabla^{\mathrm{s}} \mathbf{u}_f^{n,m} - p_f^{n,m} \boldsymbol{I} ) \cdot \mathbf{n})_\tau = \xi_f (\mathbf{u}_f^{n,m})_\tau && \mbox{at } \Gamma\,, \\[3pt]
- \mathbf{n} \cdot ( 2\mu_f \nabla^{\mathrm{s}} \mathbf{u}_f^{n,m} - p_f^{n,m} \boldsymbol{I} ) \cdot \mathbf{n} -\alpha_f \mathbf{u}_f^{n,m} \cdot \mathbf{n} =  p_p^{n,m-1} + \alpha_f \,\boldsymbol{\eta}_p \nabla p_p^{n,m-1} \cdot \mathbf{n} && \mbox{at } \Gamma\,,
\end{array}
\end{equation}
\item
the Darcy problem: find $p_p^{n,m}$ such that
\begin{equation}\label{eq:Darcy_problem_v2}
\begin{array}{l l l}
\displaystyle \frac{S_p}{\Delta t} \, p_p^{n,m} - \theta \, \nabla \cdot (\boldsymbol{\eta}_p \nabla p_p^{n,m}) = 0 && \mbox{in } \Omega_p\,, \\[5pt]
p_p^{n,m} - \alpha_p \, \boldsymbol{\eta}_p \nabla p_p^{n,m} \cdot \mathbf{n} = - \mathbf{n} \cdot ( 2\mu_f \nabla^{\mathrm{s}} \mathbf{u}_f^{n,m} - p_f^{n,m} \boldsymbol{I} ) \cdot \mathbf{n} + \alpha_p \, \mathbf{u}_f^{n,m} \cdot \mathbf{n} && \mbox{at } \Gamma\,.
\end{array}
\end{equation}

\end{enumerate}

We introduce the setting considered, e.g., in Sect. 3.2 of \cite{Discacciati:2018:IMAJNA}. Let $\Omega_f$ be the half plane $\Omega_f = \{ (x,y)\in\mathbb{R}^2 \, : \, x<0\}$, $\Omega_p$ be the complementary half plane $\Omega_p = \{ (x,y)\in\mathbb{R}^2 \, : \, x>0\}$, and $\Gamma = \{ (x,y)\in\mathbb{R}^2 \, : \, x=0\}$, so that $\mathbf{n} = (1,0)$, $\boldsymbol{\tau} = (0,1)$. Finally, let $\mu_f$ be constant in $\Omega_f$, $\boldsymbol{\eta}_p = \mbox{diag}(\eta_1,\eta_2)$ with $\eta_1$ and $\eta_2$ constant, and let $\mathbf{u}_f(x,y) = (u_1(x,y),u_2(x,y))^T$.

Using this setting and omitting the upper index $n$ for simplicity of notation, the algorithm \eqref{eq:Stokes_problem_v2}--\eqref{eq:Darcy_problem_v2} becomes: for $m\geq 1$ until convergence, solve
\begin{enumerate}
\item 
the Stokes problem
\begin{eqnarray}
\frac{1}{\Delta t}
\left( \begin{array}{l}
u_1^{m} \\
u_2^{m} 
\end{array} \right)
- \theta \, \mu_f \left( \begin{array}{l}
 \left(\partial_{xx} + \partial_{yy} \right) u_1^{m} \\
 \left( \partial_{xx} + \partial_{yy}\right) u_2^{m}
\end{array} \right)
+ \left( \begin{array}{l}
\partial_x p_f^{m} \\
\partial_y p_f^{m} 
\end{array} \right) 
= \mathbf{0} && \mbox{in }  (-\infty, 0)\times \mathbb{R}\,,  \label{eq:stokes1momentum}\\
\partial_x u_1^{m} + \partial_y u_2^{m} = 0  && \mbox{in }  (-\infty, 0)\times \mathbb{R}\,, \label{eq:stokes2continuity} \\
- \mu_f (\partial_x u_2 + \partial_y u_1) = \xi_f \,u_2^{m} && \mbox{at } \{0\}\times \mathbb{R}\,, \label{eq:stokes3bjs} \\
(-2 \mu_f \,\partial_x u_1^{m} + p_f^{m}) - \alpha_f \,u_1^{m} = p_p^{m-1} - \alpha_f (-\eta_1\, \partial_x p_p^{m-1}) && \mbox{at }  \{0\}\times \mathbb{R}\,, \label{eq:ic1rr}
\end{eqnarray}

\item
the Darcy problem
\begin{eqnarray}
\frac{S_p}{\Delta t}\, p_p^m
- \theta \, (\partial_x (\eta_1 \,\partial_x) + \partial_y (\eta_2 \,\partial_y))\, p_p^{m} =0  &&\mbox{in } (0,\infty)\times \mathbb{R}\,, \label{eq:darcy1} \\ 
p_p^m + \alpha_p ( - \eta_1\, \partial_x p_p^m) = (-2 \mu_f \,\partial_x u_1^m + p_f^m) + \alpha_p \,u_1^m && \mbox{at }  \{0\}\times \mathbb{R}\,. \label{eq:ic2rr}
\end{eqnarray}

\end{enumerate}

For the convergence analysis, similarly to \cite{Discacciati:2018:IMAJNA}, we apply the Fourier transform in the direction tangential to the interface. In this case, for ${w}(x,y) \in L^2(\mathbb{R}^2)$, the Fourier transform is defined as
\begin{equation*}
{w}(x,y) \mapsto \widehat{{w}}(x,k) = \int_{\mathbb{R}} e^{-iky} \, {w}(x,y) \, dy\,,
\end{equation*}
where $k$ is the frequency variable.

\smallskip

The following results holds.
\begin{proposition1}\label{prop:convergenceRate}
The convergence factor of the iterative method \eqref{eq:stokes1momentum}--\eqref{eq:ic2rr} does not depend on the iteration $m$ and it can be characterized for every $k\neq 0$ as $|\rho(k,\alpha_f,\alpha_p)|$ with
\begin{equation}\label{eq:convergenceFactorRho}
\rho(k,\alpha_f,\alpha_p) = 
\frac{\mathcal{G}(k)-\alpha_f}{\mathcal{G}(k)+\alpha_p} \cdot
\frac{\mathcal{H}(k)-\alpha_p}{\mathcal{H}(k)+\alpha_f}
\end{equation}
where
\begin{eqnarray}
\mathcal{G}(k) &=& \eta_p^{-1} \left( k^2 + \frac{S_p}{\eta_2\, \theta\, \Delta t} \right)^{-\frac{1}{2}} \, , \label{eq:factorG}\\
\mathcal{H}(k) &=& \frac{\displaystyle 2\mu_f \left( k^2 + \frac{1}{\mu_f \, \theta \, \Delta t} \right)^{\frac{1}{2}} - ( 2\mu_f k^2 \, \Delta t + 1) \, \mathcal{F}(k)}{1 - |k|\, \Delta t \, \mathcal{F}(k)} \, , \label{eq:factorH}
\end{eqnarray}
and
\begin{equation}\label{eq:factorF}
\mathcal{F}(k) = \frac{\displaystyle \mu_f k^2 + \left( \mu_f \left( k^2 + \frac{1}{\mu_f \, \theta \, \Delta t} \right)^{\frac{1}{2}} + \xi_f \right) \, \left( k^2 + \frac{1}{\mu_f \, \theta \, \Delta t} \right)^{\frac{1}{2}}}{k^2 \Delta t \, (2 \mu_f |k| + \xi_f)} \, ,
\end{equation}
where $\eta_p = \sqrt{\eta_1\eta_2}$.
\end{proposition1}

\medskip

\emph{Proof.} Following analogous steps as in the proof of Proposition 3.1 of \cite{Discacciati:2018:IMAJNA}, the Stokes pressure in the frequency space can be written as
\begin{equation}\label{eq:hatpf}
\hat{p}_f^{m}(x,k) = P^m(k) e^{|k|x} \, .
\end{equation}

The Fourier transform of the first component of the Stokes momentum equation \eqref{eq:stokes1momentum} becomes
\begin{equation*}
\frac{\hat{u}_1^m}{\Delta t}
- \theta \, \mu_f \partial_{xx} \hat{u}_1^{m}
+ \theta \, \mu_f k^2 \hat{u}_1^m
+ \partial_x \hat{p}_f^{m} 
= 0 \, ,
\end{equation*}
or, equivalently,
\begin{equation*}
\partial_{xx} \hat{u}_1^{m}
- \left( k^2 + \frac{1}{\mu_f \, \theta \, \Delta t} \right) \hat{u}_1^m
= \frac{|k|}{\mu_f\, \theta} P^m(k)\, e^{|k|x} \, .
\end{equation*}
The general solution of this equation is
\begin{equation}\label{eq:hatu1m}
\hat{u}_1^m(x,k) = A^m(k) \, e^{(k^2+\frac{1}{\mu_f\,\theta\,\Delta t})^{1/2} x} - |k|\, \Delta t \, P^m(k) \, e^{|k|x} \, .
\end{equation}

The Fourier transform of the second component of the Stokes momentum equation \eqref{eq:stokes1momentum} becomes
\begin{equation*}
\frac{\hat{u}_2^m}{\Delta t}
- \theta \, \mu_f \partial_{xx} \hat{u}_2^{m}
+ \theta \, \mu_f k^2 \hat{u}_2^m
+ ik \,\hat{p}_f^{m} 
= 0 \, ,
\end{equation*}
or, equivalently,
\begin{equation*}
\partial_{xx} \hat{u}_2^{m}
- \left( k^2 + \frac{1}{\mu_f \, \theta \, \Delta t} \right) \hat{u}_2^m
= \frac{ik}{\mu_f \, \theta} \, P^m(k) \, e^{|k| x} \, .
\end{equation*}
The general solution of this equation is
\begin{equation}\label{eq:hatu2m}
\hat{u}_2^m(x,k) = B^m(k) e^{(k^2+\frac{1}{\mu_f\, \theta \, \Delta t})^{1/2} x} - ik\,\Delta t \, P^m(k) \, e^{|k| x} \, .
\end{equation}

Substituting \eqref{eq:hatu1m} and \eqref{eq:hatu2m} into the Fourier transform of the Stokes continuity equation \eqref{eq:stokes2continuity}:
\begin{equation*}
\partial_x \hat{u}_1^m + ik \hat{u}_2^m = 0 \, ,
\end{equation*}
we can find the relationship
\begin{equation}\label{eq:AB}
B^m(k) = \frac{i}{k} \, A^m(k) \left( k^2 + \frac{1}{\mu_f\, \theta \, \Delta t} \right)^{1/2} \, .
\end{equation}

Moreover, the Fourier transform of the Beavers-Joseph-Saffman condition \eqref{eq:stokes3bjs} is
\begin{equation*}
-\mu_f ( \partial_x \hat{u}_2^m + ik\, \hat{u}_1^m) = \xi_f \, \hat{u}_2^m \, .
\end{equation*}
Using \eqref{eq:hatu1m} and \eqref{eq:hatu2m} at $\{0\} \times \mathbb{R}$, and the relationship \eqref{eq:AB}, we find
\begin{equation}\label{eq:PAF}
P^m(k) = A^m(k) \, \mathcal{F}(k)
\end{equation}
with $\mathcal{F}(k)$ defined as in \eqref{eq:factorF}.

\smallskip

Taking now the Fourier transform of the Darcy equation \eqref{eq:darcy1}, we find
\begin{equation*}
\frac{S_p}{\Delta t} \, \hat{p}_p^{m} - \theta \, \eta_1 \partial_{xx} \hat{p}_p^m + \theta \, \eta_2 \, k^2 \hat{p}_p^{m} = 0 \, ,
\end{equation*}
or, equivalently,
\begin{equation*}
\partial_{xx} \hat{p}_p^m - \frac{\eta_2}{\eta_1} \left( k^2 + \frac{S_p}{\eta_2 \, \theta \, \Delta t} \right) \hat{p}_p^m = 0 \, ,
\end{equation*}
whose general solution is
\begin{equation}\label{eq:hatpp}
\hat{p}_p^m(x,k) = \Phi^m(k) e^{- (\frac{\eta_2}{\eta_1} ( k^2 + \frac{S_p}{\eta_2 \, \theta \, \Delta t}))^{1/2} x} \, .
\end{equation}

The Fourier transforms of the interface conditions \eqref{eq:ic1rr} and \eqref{eq:ic2rr} become
\begin{eqnarray}
&& (-2 \mu_f \partial_x \hat{u}_1^{m} + \hat{p}_f^{m}) - \alpha_f \,\hat{u}_1^{m} = \hat{p}_p^{m-1} - \alpha_f (-\eta_1\, \partial_x \hat{p}_p^{m-1}) \, , \label{eq:ic1rr_fourier}\\
&& \hat{p}_p^m + \alpha_p ( - \eta_1\, \partial_x \hat{p}_p^m) = (-2 \mu_f \partial_x \hat{u}_1^m + \hat{p}_f^m) + \alpha_p \,\hat{u}_1^m \, . \label{eq:ic2rr_fourier}
\end{eqnarray}

Using the expressions \eqref{eq:hatpf}--\eqref{eq:hatu2m}, \eqref{eq:hatpp} at $\{0\}\times\mathbb{R}$ and the relationships \eqref{eq:AB} and \eqref{eq:PAF}, from \eqref{eq:ic1rr_fourier} and \eqref{eq:ic2rr_fourier} we obtain, respectively,
\begin{eqnarray}
&& \left( - \left( \alpha_f + 2\mu_f \left( k^2 + \frac{1}{\mu_f\,\theta \, \Delta t} \right)^{\frac{1}{2}} \right) + (2\mu_f k^2 \Delta t + \alpha_f |k| \Delta t +  1 ) \, \mathcal{F}(k) \right) A^m(k) \nonumber \\
&& \qquad = \left( 1 - \alpha_f \, \eta_p \left( k^2 + \frac{S_p}{\eta_2\, \theta \, \Delta t} \right)^{\frac{1}{2}} \right) \Phi^{m-1}(k)\,, \label{eq:itermet1bis} \\
&& \left( 1 + \alpha_p \,\eta_p \left( k^2 + \frac{S_p}{\eta_2\, \theta\, \Delta t} \right)^{\frac{1}{2}} \right) \Phi^{m}(k) \nonumber \\
&& \qquad = \left( \left( \alpha_p - 2\mu_f \left( k^2 + \frac{1}{\mu_f\,\theta\,\Delta t} \right)^{\frac{1}{2}} \right) + (2\mu_f k^2 \Delta t - \alpha_p |k| \Delta t + 1 ) \, \mathcal{F}(k) \right) A^m(k) \,.\label{eq:itermet2bis}
\end{eqnarray}

\smallskip

The convergence factor \eqref{eq:convergenceFactorRho} can now be directly obtained from \eqref{eq:itermet1bis} and \eqref{eq:itermet2bis}, using the definitions \eqref{eq:factorG} and \eqref{eq:factorH}.

\smallskip

A direct calculation shows that $1 - |k|\, \Delta t \, \mathcal{F}(k) < 0$ for all $k \neq 0$, so that $\mathcal{H}(k)$ is well defined.
\hfill $\Box$

\medskip

%

\begin{remark}
Proposition \ref{prop:convergenceRate} explicitly excludes the frequency $k=0$. In fact, in this case, the only admissible Stokes pressure $\hat{p}_f^m(x,k)$ would be the zero function since it must satisfy the homogeneous Laplace equation $-\partial_{xx}\hat{p}_f^m=0$ in the infinite half plane $\Omega_f$, and decay at infinity in order to be $L^2$ integrable.
The velocities would become $\widehat{u}_1^m(x,k)=A^m(k)e^{\sqrt{\frac{1}{\mu_f\theta\Delta t}} x}$, and $\widehat{u}_2^m(x,k)=B^m(k)e^{\sqrt{\frac{1}{\mu_f\theta\Delta t}} x}$. Inserting these expressions into the divergence-free constraint and into the Beavers-Joseph-Saffman condition finally leads to $\widehat{u}_1^m(x,k)=\widehat{u}_2^m(x,k)=0$. 
The frequency $k=0$ must be accounted for if one considers periodic boundary conditions in a bounded domain in the $y$-direction as in \cite{Gander:2020}. In such case, due to boundedness of the domain, one obtains $\widehat{u}^m_1=\widehat{u}^m_2=0$, $\widehat{p}_f^m(x,k)=P^m\in \mathbb{R}$, and $\widehat{p}_p^m(x,k)=\Phi^m(k) e^{-\sqrt{\frac{S_p}{\theta\eta_1 \Delta t}}x}$. Inserting these expressions into the transmission conditions, we obtain $\rho(k=0,\alpha_f,\alpha_p)=\frac{1-\alpha_f\eta_p \sqrt{\frac{S_p}{{\theta\eta_2 \Delta t}}}}{1+\alpha_p\eta_p \sqrt{\frac{S_p}{{\theta\eta_2 \Delta t}}}}.$
\end{remark}

Now, our goal is to characterize optimal parameters $\alpha_f$ and $\alpha_p$ that minimize the convergence factor \eqref{eq:convergenceFactorRho}. To this aim, we study the min-max problem
\begin{equation}\label{eq:minmax_problem1}
\min_{\alpha_f,\alpha_p\in \mathbb{R}^+}\max_{[k_{\min},k_{\max}]} |\rho(k,\alpha_f,\alpha_p)| \, ,
\end{equation}
where $[k_{\min},k_{\max}]$ is the range of frequencies of interest for the problem, and they are usually approximated by $k_{\min}=\frac{\pi}{|\Gamma|}$ and $k_{\max}=\frac{\pi}{h}$, where $h$ is the mesh size \cite{Gander:2006:OSM}.

\subsection{Optimization of the convergence factor in a simplified setting}\label{sect:Optimization}

The solution of the min-max problem \eqref{eq:minmax_problem1} is quite challenging, and it significantly differs from several min-max problems studied in the literature (see \cite{Gander:2006:OSM,Gander:twosided,Gander_XU:2016:OSM,gander2019heterogeneous} and references therein). Moreover, the expression of $\mathcal{F}(k)$ is quite complex and it prevents the use of direct calculations. Therefore, in order to make the min-max problem \eqref{eq:minmax_problem1} feasible of a theoretical study, we replace the Beaver-Joseph-Saffman condition \eqref{eq:interf3} by the zero tangential interface velocity 
\begin{equation}\label{eq:interf3Simplified}
(\mathbf{u}_f)_\tau = 0 \quad \mbox{on } \Gamma \, ,
\end{equation}
(see, e.g., in \cite{Dangelo:2011:M2AN}), which corresponds to the limit $\xi_f\rightarrow \infty$ in \eqref{eq:interf3}. We find
\begin{equation*}
\lim_{\xi_f\rightarrow \infty} \mathcal{F}(k)=\frac{\displaystyle \left( k^2+\frac{1}{\mu_f \, \theta \, \Delta t} \right)^{\frac{1}{2}}}{k^2 \, \Delta t} \, .
\end{equation*}
Using this limit in \eqref{eq:factorH}, we obtain the simplified expression
\begin{equation}\label{eq:hSimplified}
\mathcal{H}(k) =
\frac{\displaystyle \left( k^2+\frac{1}{\mu_f \, \theta \, \Delta t} \right)^{\frac{1}{2}}}{\displaystyle |k| \, \Delta t \, \left( \left( k^2+\frac{1}{\mu_f \, \theta \, \Delta t} \right)^{\frac{1}{2}} - |k| \right)} \, .
\end{equation}


\medskip

We now aim to find the optimal $\alpha_f$ and $\alpha_p$ that solve the min-max problem \eqref{eq:minmax_problem1} under the simplifying assumption \eqref{eq:hSimplified}.

In \cite{gander2019heterogeneous,gander2017optimized}, the authors propose to obtain the solution of \eqref{eq:minmax_problem1} with $\rho(k,\alpha_f,\alpha_p)$ in the form \eqref{eq:convergenceFactorRho} by setting $\alpha_f=\mathcal{G}(s)$ and $\alpha_p=\mathcal{H}(s)$, for a $s\in \mathbb{R}^+$. If the functions $\mathcal{G}(\cdot)$ and $\mathcal{H}(\cdot)$ are strictly monotonic, the convergence factor has only one zero for $k>0$ located at $k=s$. Then, analyzing the derivative of $\rho(k,\alpha_f(s),\alpha_p(s))$ with respect to $s$, one can conclude that the optimal parameter, say $s_{opt}$, is the one leading to an equioscillation between the convergence factor at $k_{\min}$ and $k_{\max}$, i.e., $|\rho(k_{\min},\alpha_f(s_{opt}),\alpha_p(s_{opt}))| = |\rho(k_{\max},\alpha_f(s_{opt}),\alpha_p(s_{opt}))|$. (This argument can be extended to double-sided optimized transmission conditions, letting $\alpha_f=\mathcal{G}(s)$ and $\alpha_p=\mathcal{H}(p)$ for some $s,p\in \mathbb{R}^+$.) The key point of this strategy is that both $\mathcal{G}(\cdot)$ and $\mathcal{H}(\cdot)$ must be strictly monotonic, so that the convergence factor has as many zeros as the parameters to optimize.

\smallskip

For the problem at hand, this approach is not viable. Indeed, the function $k\mapsto \mathcal{G}(k)$ in \eqref{eq:factorG} is strictly positive and decreasing $\forall k>0$, but the map $k\mapsto \mathcal{H}(k)$ in \eqref{eq:hSimplified} is strictly positive, decreasing for $k<\widehat{k}$ and increasing for $k> \widehat{k}$ where
\begin{equation}\label{eq:kHat}
\widehat{k}  = \frac{\sqrt{2}}{2}\sqrt{\frac{\sqrt{5}-1}{\mu_f \, \theta \, \Delta t}} \, ,
\end{equation}
with $\lim_{k\rightarrow +\infty} \mathcal(k)=+\infty$ and $\lim_{k\rightarrow 0^+} \mathcal(k)=+\infty$. 
Hence, setting $\alpha_f=\mathcal{G}(s)$ and $\alpha_p=\mathcal{H}(s)$ would lead to two zeros for $k>0$, one at $k=s$ (for which both $\mathcal{G}(k)-\mathcal{G}(s)=0$ and $\mathcal{H}(k)-\mathcal{H}(s)=0$) and another one at a point, say $k^{\star}(s)$, which is the other preimage of $\mathcal{H}(s)$, i.e., the unique other point satisfying $\mathcal{H}(k^{\star}(s)) = \mathcal{H}(s)$. Only for $s=\widehat{k}$, one has $k^{\star}(\widehat{k})=\widehat{k}$.

The additional difficulty of performing symbolic calculations with the expressions of $\mathcal{G}$ and $\mathcal{H}$ suggests us to simplify the problem by setting $\alpha_p = \mathcal{H}(\widehat{k})=:l$, that is we set $\alpha_p$ equal to the value of $\mathcal{H}(\cdot)$ at the minimum $\widehat{k}$. In this way, the convergence factor $\rho(k,s)$ has now two zeros, one fixed at $k=\widehat{k}$, the other at $k=s$, and thus only one zero is parameter-dependent.

\begin{remark}
Setting $\alpha_p=l$ does not necessarily lead to less efficient transmission conditions than choosing $\alpha_p=\mathcal{H}(s)$ and then optimizing also $\alpha_p$. In fact, $\rho(k,s)$ has now one zero which is fixed and it depends on the physical and discretization parameters $\mu_f$, $\theta$ and $\Delta t$, and the optimization of $\alpha_f$ automatically takes into account the position of this fixed zero. In the case when both $\alpha_f$ and $\alpha_p$ depend on the the free parameter $s$, one apparently optimizes both $\alpha_p$ and $\alpha_f$, but actually these parameters are constrained to move on the one dimensional curve $\Sigma := \left\{ (x,y)\in \mathbb{R}^+ \,:\, (x,y) = (\mathcal{G}(s),\mathcal{H}(s)), \; s\in \mathbb{R}^+\right \},$ see, e.g., \cite{Discacciati:2018:IMAJNA}.
\end{remark}

Thus, we consider the min-max problem
\begin{equation}\label{eq:min_max}
\min_{s\in \mathbb{R}^+} \max_{k\in [k_{\min},k_{\max}]} \rho(k,s) = \min_{s\in \mathbb{R}^+} \max_{k\in [k_{\min},k_{\max}]} 
\left|\frac{\mathcal{G}(k)-\mathcal{G}(s)}{\mathcal{H}(k)+\mathcal{G}(s)} \cdot \frac{\mathcal{H}(k)-l}{\mathcal{G}(k)+l}\right|,
\end{equation}
where $l=\mathcal{H}(\widehat{k})$ and $\widehat{k}$ is defined in \eqref{eq:kHat}.
The following Theorem characterizes the optimal solution of \eqref{eq:min_max} according to the position of $\widehat{k}$ with respect to $[k_{\min},k_{\max}]$. 

\begin{theorem}\label{thm:minmax}
\begin{itemize}
\item Assume that $\widehat{k}\in  [k_{\min},k_{\max}]$.
The solution of \eqref{eq:min_max} is given by
\begin{itemize}
\item The unique $s_1^{\star}$ such that $|\rho(k_{\min},s^{\star}_1)|=|\rho(k_{\max},s^{\star}_1)|=:M$, if
$|\rho(\widetilde{k}(s_1^{\star}),s_1^{\star})|\leq M$, where $\widetilde{k}$ is the unique local interior maximum.
\item If $|\rho(\widetilde{k}(s_1^{\star}),s_1^{\star})|>M$ and $s_1^{\star}< \widehat{k}$, then the optimal solution is given by $s_2^{\star}$ which is the unique solution of $|\rho(k_{\min},s)|=|\rho(\widetilde{k},s)|$ for $s_1^{\star}<s<\widehat{k}$.
\item If $|\rho(\widetilde{k}(s_1^{\star}),s_1^{\star})|>M$ and $s_1^{\star}> \widehat{k}$,
then the optimal solution is given by $s_3^{\star}$ which is the unique solution of $|\rho(k_{\max},s)|=|\rho(\widehat{k},s)|$ for $\widehat{k}<s<s_1^{\star}$. 
\end{itemize}
\item Assume that $\widehat{k}> k_{\max}$.  The solution of \eqref{eq:min_max} is given by
\begin{itemize}
\item The unique $s_1^{\star}$ such that $|\rho(k_{\min},s^{\star}_1)|=|\rho(k_{\max},s^{\star}_1)|$, if $\widetilde{k}(s_1^{\star})>k_{\max}$.
\item If $\widetilde{k}(s_1^{\star})<k_{\max}$, then the optmimal solution is given by
$s_4^{\star}$ which is the unique solution of $|\rho(k_{\min},s)|=|\rho(\widetilde{k},s)|$.
\end{itemize}
\item Assume that $\widehat{k}< k_{\min}$.  The solution of \eqref{eq:min_max} is given by
\begin{itemize}
\item The unique $s_1^{\star}$ such that $|\rho(k_{\min},s^{\star}_1)|=|\rho(k_{\max},s^{\star}_1)|$, if $\widetilde{k}(s_1^{\star})<k_{\min}$.
\item If $\widetilde{k}(s_1^{\star})>k_{\min}$, then the optmimal solution is given by
$s_5^{\star}$ which is the unique solution of $|\rho(\widetilde{k},s)|=|\rho(k_{\max},s)|$.
\end{itemize}
\end{itemize}
\end{theorem}

\begin {figure}
\begin{tikzpicture}[scale=0.51]
\draw[->,ultra thick] (0,0)--(6,0) node[right]{$k$};
\draw[->,ultra thick] (0,0)--(0,6) node[left]{$|\rho(k,s)|$};
\draw[dashed] (1,0)--(1,5); 
\node at (1,5.3) {$k_{\min}$};
\draw[dashed] (5,0)--(5,5); 
\node at (5,5.3) {$k_{\max}$};
\node at (1.7,-0.8) {$\widehat{k}$};
\draw[black,ultra thick] (1.7,-0.2)--(1.7,0.2); 
\draw[black,ultra thick] (3,-0.2)--(3,0.2); 
\node at (3,-1) {$s_1^\star$};
\draw[thin,black] plot [smooth] coordinates { (0.2,4)  (1,1.6) (1.7,0) (2.4,1) (3,0) (5,1.6) (5.5,2)};
\draw[dotted,thick] (1,1.6)--(5,1.6); 
\draw [->,thick,black] (1.2,1.6) -- (1.2,2.2);
\draw [->,thick,black] (2.4,1.6) -- (2.4,2.2);
\draw [<-,thick,black] (4.8,1.6) -- (4.8,2.2);
\end{tikzpicture}
\begin{tikzpicture}[scale=0.51]
\draw[->,ultra thick] (0,0)--(6,0) node[right]{$k$};
\draw[->,ultra thick] (0,0)--(0,6) node[left]{$|\rho(k,s)|$};
\draw[dashed] (1,0)--(1,5); 
\node at (1,5.3) {$k_{\min}$};
\draw[dashed] (5,0)--(5,5); 
\node at (5,5.3) {$k_{\max}$};
\node at (1.7,-0.8) {$\widehat{k}$};
\draw[black,ultra thick] (1.7,-0.2)--(1.7,0.2); 
\draw[black,ultra thick] (3,-0.2)--(3,0.2); 
\node at (3,-1) {$s_1^\star$};
\draw[thin,black] plot [smooth] coordinates { (0.2,4)  (1,1.6) (1.7,0) (2.4,2) (3,0) (5,1.6) (5.5,2)};
\draw[dotted,thick] (1,1.6)--(5,1.6); 
\draw [->,thick,black] (1.2,1.6) -- (1.2,2.2);
\draw [->,thick,black] (2.4,2.2) -- (2.4,2.8);
\draw [<-,thick,black] (4.8,1.6) -- (4.8,2.2);
\end{tikzpicture}
\begin{tikzpicture}[scale=0.51]
\draw[->,ultra thick] (0,0)--(6,0) node[right]{$k$};
\draw[->,ultra thick] (0,0)--(0,6) node[left]{$|\rho(k,s)|$};
\draw[dashed] (1,0)--(1,5); 
\node at (1,5.3) {$k_{\min}$};
\draw[dashed] (5,0)--(5,5); 
\node at (5,5.3) {$k_{\max}$};
\node at (1.7,-0.8) {$\widehat{k}$};
\draw[black,ultra thick] (1.7,-0.2)--(1.7,0.2); 
\draw[black,ultra thick] (2.5,-0.2)--(2.5,0.2); 
\node at (2.5,-1) {$s_3^\star$};
\draw[thin,black] plot [smooth] coordinates { (0.2,4)  (1,0.8) (1.7,0) (2,1.4) (2.5,0) (5,1.4) (5.5,2)};
\draw[dotted,thick] (1,1.4)--(5,1.4); 
\end{tikzpicture}
\caption{The three panels support the proof of Proposition \ref{thm:minmax}. The black arrows show how $|\rho(k,s)|$ behaves in the three local maxima as $s$ increases. In the left panel, equioscillating between $|\rho(k_{\min},s)|$ and $|\rho(k_{\max},s)|$ leads to an optimal solution. In the central panel, $s_1^\star$ is not optimal, as decreasing $s$ would decrease $|\rho|$ at both $k_{\min}$ and $\widetilde{k}$ and increase it at $k_{\max}$ until $|\rho(\widetilde{k},s)|=|\rho(k_{\max},s)|$}\label{Fig:thm}
\end{figure}
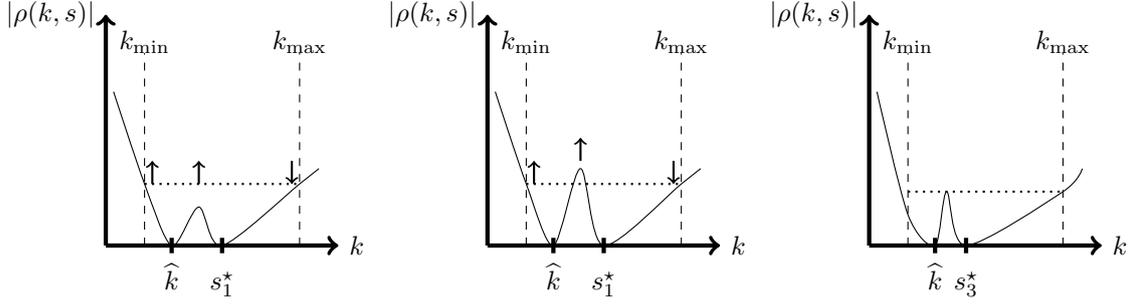

\begin{proof}
Let us start with a preliminary analysis of $\mathcal{H}(\cdot)$ and $\rho(\cdot,\cdot)$.
As observed above, $\mathcal{H}(k)$ is a strictly decreasing function for $k<\widehat{k}$, strictly increasing for $k>\widehat{k}$, and with minimum at $k=\widehat{k}$. Thus, $(\mathcal{H}(k)-l)$ is nonnegative for every $k$, and $\rho(k,s)$ has two zeros, the one first at $k=s$ and the second one fixed at $k=\widehat{k}$. Further, $\rho(k,s)$ is positive for $k< s$ and negative for $k>s$.

The sign of the partial derivative with respect to $s$ is equal to 
\[
  \sign\left(\frac{\partial |\rho(k,s)|}{\partial s}\right) 
  =
  \sign(\rho(k,s)) \;
  \sign\left(\frac{\mathcal{G}^\prime (s)(l-\mathcal{H}(k))(\mathcal{H}(k)+\mathcal{G}(k))}{(\mathcal{G}(k)+l)(\mathcal{H}(k)+\mathcal{G}(s))^2}\right).
\]

Being the second term on the right hand side always positive since $(l-\mathcal{H}(k))<0$ and $\mathcal{G}$ is strictly decreasing, there holds
\[
  \sign\left(\frac{\partial |\rho(k,s)|}{\partial s} \right) = 
  \sign(\rho(k,s)) = \sign(s-k) \, .
\]
We conclude that the optimal $s$ must be greater than $k_{\min}$, since otherwise we would have $\sign\left(\frac{\partial |\rho(k,s)|}{\partial s}\right)<0$ for every $k\in [k_{\min},k_{\max}]$, hence we could not be at the optimum, since increasing $s$ would decrease $|\rho(k,s)|$ for all $k\in [k_{\min},k_{\max}]$. For analogous reasons, the optimal $s$ must be less than $k_{\max}$, so that we can restrict the range of $s$ to the interval $[k_{\min},k_{\max}]$. 
Studying the derivative with respect to $k$, one notices that $|\rho(k,s)|$ is decreasing for $k\leq \min\{s,\widehat{k}\}$ and increasing for
$k\geq \max \{s,\widehat{k}\}$. Between the two zeros, it has a unique local maximum in $\widetilde{k}(s)$.

Let us now suppose that $\widehat{k}\in [k_{\min},k_{\max}]$. Since $s\in [k_{\min},k_{\max}]$, the local maximum $\widetilde{k}(s)$ surely lies in $[k_{\min},k_{\max}]$. The min-max problem \eqref{eq:min_max} then simplifies to
\[\min_{s\in \mathbb{R}^+} \max_{k\in [k_{\min},k
_{\max}]} |\rho(k,s)|= \min_{s\in [k_{\min},k_{\max}]} \max\left\{|\rho(k_{\min},s)|,\rho(\widetilde{k}(s),s)|,|\rho(k_{\max},s)|\right\}.\]
We now study how $|\rho(k_{\min},s)|$ and  $|\rho(k_{\max},s)|$ depend on $s$. For all $s\in (k_{\min},k_{\max}]$, there holds
\[\sign\left(\frac{\partial |\rho(k,s)|}{\partial s}|_{k=k_{\min}}\right)=\sign(s-k_{\min})>0,\quad \sign\left(\frac{\partial |\rho(k,s)|}{\partial s}|_{k=k_{\max}}\right)=\sign(s-k_{\max})<0.\] 
Thus, $|\rho(k_{\min},s)|$ is a strictly increasing function of $s$ that satisfies $|\rho(k_{\min},k_{\min})|=0$ and $|\rho(k_{\min},k_{\max})|>0$. Similarly, $|\rho(k_{\max},s)|$ is strictly decreasing, with $|\rho(k_{\max},k_{\min})|>0$ and $|\rho(k_{\max},k_{\max})|=0$. By continuity there exists a unique $s_1^{\star}$ such that $|\rho(k_{\min},s_1^{\star})|=|\rho(k_{\max},s_1^{\star})|=:M$. Now, if $|\rho(\widetilde{k}(s_1^{\star}),s_1^{\star})|\leq M$, we have found the optimum, since perturbing $s$ would increase either $|\rho(k_{\min},s)|$, or $|\rho(k_{\max},s)|$, and thus the maximum of $|\rho(k,s)|$, see left panel of Fig. \ref{Fig:thm}. Suppose instead that $|\rho(\widetilde{k}(s_1^{\star}),s_1^{\star})|> M$. In this case, the optimal solution is obtained by equioscillating $|\rho(\widetilde{k}(s),s)|$ either with $|\rho(k_{\min},s)|$ or $|\rho(k_{\max},s)|$, since they are respectively strictly increasing and strictly decreasing with respect to $s$ (central panel of Fig. \ref{Fig:thm}). Further, the total derivative of $|\rho(\widetilde{k}(s),s)|$ with respect to $s$ is
\[\frac{d|\rho(\widetilde{k}(s),s)|}{ds}=\sign(s-k)\left(\frac{\partial |\rho(k,s)|}{\partial k}|_{k=\widetilde{k}(s)}\cdot \widetilde{k}^\prime(s) + \frac{\partial \rho(\widetilde{k}(s),s)}{\partial s}\right)=\sign(s-k)\frac{\partial \rho(\widetilde{k}(s),s)}{\partial s},\]
whose sign is
\[\sign\left(\frac{d|\rho(\widetilde{k}(s),s)|}{ds}\right)=\sign(s-k).\]
If $s_1^{\star}<\widehat{k}$ and $|\rho(\widetilde{k}(s_1^{\star}),s_1^{\star})|> M$,
decreasing $s$ in the range $[k_{\min},s_1^{\star}]$ would not improve the maximum of the convergence factor as $|\rho(\widetilde{k}(s),s)|$ increases. On the other hand, increasing $s$ in the range $[s_1^{\star},\widehat{k}]$ would increase $|\rho(k_{\min},s)|$ and decrease both $|\rho(\widetilde{k},s)|$ and $|\rho(k_{\max},s)|$. There exists surely a $s_2^{\star} \in [s_1^{\star},\widehat{k}]$ for which $|\rho(\widetilde{k}(s_2^{\star}),s_2^{\star})|=|\rho(k_{\min},s_2^{\star})|$ since  $|\rho(\widetilde{k}(\widehat{k}),\widehat{k})|=0$. The optimal solution is then $s_2^{\star}$ as increasing further $s$ would increase $|\rho(k_{\min},s)|$ as well. Similarly, one shows that if $s_1^{\star}>\widehat{k}$, the optimal solution is given by $s_3^{\star}\in [\widehat{k},s_1^{\star}]$ such that $|\rho(\widetilde{k}(s),s)|=|\rho(k_{\max},s)|$ (right panel of Fig. \ref{Fig:thm}). This ends the proof of the first claim.

Next, suppose that $\widehat{k}> k_{\max}$, so that the local maximum $\widetilde{k}(s)$ is not necessarily in $[k_{\min},k_{\max}]$. Let us define $s_4^{\star}\in [k_{\min},k_{\max}]$ be the unique solution of $|\rho(k_{\min},s)|=|\rho(k_{\max},s)|$. If $\widetilde{k}(s_4^{\star})>k_{\max}$, then $s_4^{\star}$ is the optimal solution as varying $s$ would increase either $|\rho(k_{\min},s)|$ or $|\rho(k_{\max},s)|$. In constrast, if $\widetilde{k}(s_4^{\star})<k_{\max}$, then $|\rho(\widetilde{k}(s_4^{\star}),s_4^{\star})|>|\rho(k_{\min},s_4^{\star})|$ and it is convenient to increase $s$ until $|\rho(\widetilde{k}(s),s)|=|\rho(k_{\min},s)|$. The case $\widehat{k}<k_{\min}$ is treated similarly.

\end{proof}

\begin{remark} 
Theorem \ref{thm:minmax} considers three cases depending on the relative position of $\widehat{k}$ with respect to $[k_{\min},k_{\max}]$. This
can be easily determined and it depends on few physical, geometrical and discretization parameters. Indeed, using the standard ansatzes $k_{\min}=\frac{\pi}{|\Gamma|}$ and $k_{\max}=\frac{\pi}{h}$ \cite{Gander:2006:OSM}, from \eqref{eq:kHat} we get
\begin{equation*}
\begin{array}{ll}
\widehat{k}<k_{\min} &  \text{if}\quad \mu_f\theta\Delta t >\frac{(\sqrt{5}-1)|\Gamma|^2}{2\pi^2},\\[5pt]
k_{\min}\leq \widehat{k}\leq k_{\max}\quad &\text{if}\quad \frac{(\sqrt{5}-1)h^2}{2\pi^2}\leq \mu_f\theta\Delta t\leq \frac{(\sqrt{5}-1)|\Gamma|^2}{2\pi^2},\\[5pt]
\widehat{k}>k_{\max}\quad &\text{if}\quad \mu_f\theta\Delta t <\frac{(\sqrt{5}-1)h^2}{2\pi^2}.
\end{array}
\end{equation*}
Notice that $\widehat{k}<k_{\min}$ poses a tight constraint on the value of $\Delta t$, which is not likely to be satisfied in practice. The other two cases are more realistic, and the specific regime depends on the choice of $\Delta t$ and $h$ which is dictated by the finite element spaces and the time integration scheme used.
\end{remark}

\section{Numerical results}\label{sect:NumericalResults}

In this section, we discuss numerical results to test the proposed framework. Subsection \ref{sec:choice_parameters} introduces the physical relevant parameters used in the simulations.
Subsection \ref{sec:periodic} validates the Fourier analysis by considering the Stokes-Darcy coupling on a bounded domain with periodic boundary conditions and with zero tangential velocity, so to mimic the hypothesis of Section \ref{sec:Robin-Robin}. Next, in subsection \ref{sec:test_analytic} we violate the assumptions of the Fourier analysis by considering nonperiodic boundary conditions, together with the Beaver-Joseph-Saffman boundary condition. Finally in subsection \ref{sec:realistic}, we present the convergence of the Robin-Robin method for a realistic application.

\subsection{Choice of physically relevant parameters}\label{sec:choice_parameters}

For the numerical tests we consider an incompressible fluid with kinematic viscosity $\nu = 10^{-6}$ m$^2$/s (water) at Reynolds $0.1 \leq Re \leq 5$, and we choose the characteristic dimension of the domain $\Omega_f$ to be $X_f = 0.05$ m. For Darcy's problem, the dimensionless coefficient $S_p$ is obtained as $S_p = \nu^2 \frac{S_0}{g} \frac{Re^2}{X_f^2}$, where $g=9.8$ $\mbox{m}\times\mbox{s}^{-2}$ is the magnitude of the gravity acceleration, and $S_0$ is the specific storage whose value ranges, e.g., in the interval $10^{-5} \,\mbox{m}^{-1} \leq S_0 \leq 10^{-3} \,\mbox{m}^{-1}$ \cite{Moraiti:2012:JMAA}. The porous medium is characterized by constant intrinsic permeability $\mathbf{K}$ in the range $10^{-11} \,\mbox{m}^{2} \leq \mathbf{K} \leq 10^{-8} \,\mbox{m}^{2}$ \cite{Bear:1979:HOG}, so that $\boldsymbol{\eta}_p = \frac{\mathbf{K}}{X_f^2}\, Re$ is also constant with $\eta_1=\eta_2$. Using the indicated values, we identify four test cases (A)--(D) with parameters of physical relevance that we will use for our numerical experiments. The test cases are reported in Table \ref{table:values} with values of the dimensionless quantities rounded to 2 decimal places.

\begin{table}[bht]
\begin{center}
\begin{tabular}{cllll}
Test & \multicolumn{1}{c}{(A)} & \multicolumn{1}{c}{(B)} & \multicolumn{1}{c}{(C)} & \multicolumn{1}{c}{(D)} \\
\hline
$Re$                  & \multicolumn{1}{c}{$0.1$} & \multicolumn{1}{c}{$1$} & \multicolumn{1}{c}{$0.1$} & \multicolumn{1}{c}{$5$} \\
$S_p$                 & $4.08 \times 10^{-16}$ & $4.08 \times 10^{-15}$ & $4.08 \times 10^{-18}$ & $1.02 \times 10^{-14}$ \\
$\boldsymbol{\eta}_p$ & $4.00 \times 10^{-10}$ & $4.00 \times 10^{-7}$ & $4.00 \times 10^{-9}$ & $2.00 \times 10^{-7}$ \\
$\xi_f$               & $1.58 \times 10^{+5}$ & $1.58 \times 10^{+3}$ & $5.00 \times 10^{+4}$ & $1.00 \times 10^{+3}$ \\
\hline
\end{tabular}
\caption{Values of the dimensionless parameters defining four test cases (A)--(D).}\label{table:values}
\end{center}
\end{table}

\subsection{Tests with periodic boundary conditions}\label{sec:periodic}
Consider the domain $\Omega=(-1,1)^2$ decomposed into $\Omega_f=(0,1)^2$ and $\Omega_p=(0,1)\times(-1,0)$, with $\Gamma=(0,1)\times \left\{0\right\}$.
The discretization is based on a uniform mesh of squares of edge length $h_j=0.1\times 2^{1-j}$ $j=1,\dots,4$. Each square is then divided into two right triangles. We use Taylor-Hood finite elements for Stokes and continuous Lagrangian $\mathbb{P}_2$ finite elements for Darcy's pressure. For the time discretization we use the Crank-Nicolson method ($\theta=\frac{1}{2}$).
We impose homogeneous Dirichlet boundary conditions on the horizontal edges of $\Omega$ and periodic boundary conditions on the vertical edges, which permit to mimic the assumption on the unboundedness of the domain required by the Fourier analysis, see \cite{Gander:2020}. 
In this setting, we consider the error equation (i.e., $\chi_f^n=0$, $\chi_p^n=0$ in \eqref{eq:interfaceSystem_v2}) and study the convergence of the stationary iteration \eqref{eq:stationary_iteration} to the zero solution starting from a random initial guess. Theorem \ref{thm:minmax} does not cover the case under study as we need to include the zero frequency $k=0$. We thus compute the optimized value of $s$ by solving numerically the min-max problem $\min_{s\in \mathbb{R}}\max_{k\in\left\{0\right\}\cup [k_{\min},k_{\max}]} |\rho(k,s)|$ using the Nelder-Mead algorithm \cite{Lagarias:1998:CPN}.
Table \ref{table:testperiodic} reports the values of the optimal parameters $\alpha_f$, $\alpha_p$, and the number of iterations to reach a tolerance of $10^{-8}$  on the relative error of the stationary iterations for different values of $h_j$ and for timesteps $\Delta t_1=0.05$, $\Delta t_2=0.01$, $\Delta t_3=0.005$ and $\Delta t_4=0.001$.
Notice that very large values of $\alpha_f$ essentially transform the Robin-Robin algorithm into a Dirichlet-Robin algorithm, in which the Stokes subdomain receives the Darcy velocity as boundary conditions at the interface $\Gamma$.  
\begin{table}[bht]
{\footnotesize
\begin{center}
\begin{tabular}{ccllcc}
 Test & Mesh & \multicolumn{1}{c}{$\alpha_f$} & \multicolumn{1}{c}{$\alpha_p$} & iter  &\\
\hline
(A) & $h_1$ & $2.50 \times 10^{+9}$ & $1.49 \times 10^{+2}$ & $4$ &\\
    & $h_2$ & $2.50 \times 10^{+9}$ & $1.49 \times 10^{+2}$ & $4$ &\\
    & $h_3$ & $2.50 \times 10^{+9}$ & $1.49 \times 10^{+2}$ & $6$ &\\
    & $h_4$ & $2.49 \times 10^{+9}$ & $1.49 \times 10^{+2}$ & $6$ &\\
\hline
(B) & $h_1$ & $2.49 \times 10^{+6}$ & $4.71 \times 10^{+1}$ & $6$ &\\
    & $h_2$ & $2.47 \times 10^{+6}$ & $4.70 \times 10^{+1}$ & $8$ &\\
    & $h_3$ & $2.39 \times 10^{+6}$ & $4.71 \times 10^{+1}$ & $12$ &\\
    & $h_4$ & $2.02 \times 10^{+6}$ & $4.71 \times 10^{+1}$ & $20$ &\\
\hline
(C) & $h_1$ & $2.50 \times 10^{+8}$ & $1.49 \times 10^{+2}$ & $6$ &\\
    & $h_2$ & $2.50 \times 10^{+8}$ & $1.49 \times 10^{+2}$ & $6$ & \\
    & $h_3$ & $2.49 \times 10^{+8}$ & $1.49 \times 10^{+2}$ & $6$ &\\
    & $h_4$ & $2.45 \times 10^{+8}$ & $1.49 \times 10^{+2}$ & $10$ &\\
\hline
(D) & $h_1$ & $5.00 \times 10^{+6}$ & $2.11 \times 10^{+1}$ & $6$ & \\
    & $h_2$ & $5.00 \times 10^{+6}$ & $2.11 \times 10^{+1}$ & $6$ & \\
    & $h_3$ & $4.98 \times 10^{+6}$ & $2.11 \times 10^{+1}$ & $6$ &\\
    & $h_4$ & $4.91 \times 10^{+6}$ & $2.11 \times 10^{+1}$ & $10$ & \\
\hline
\end{tabular}
\begin{tabular}{ccllcc}
 Test & Mesh & \multicolumn{1}{c}{$\alpha_f$} & \multicolumn{1}{c}{$\alpha_p$} & iter  &\\
\hline
(A) & $\Delta t_1$ & $2.50 \times 10^{+9}$ & $6.66 \times 10^{+1}$ & $4$ &\\
    & $\Delta t_2$ & $2.50 \times 10^{+9}$ & $1.49 \times 10^{+2}$ & $4$ &\\
    & $\Delta t_3$ & $2.50 \times 10^{+9}$ & $2.11 \times 10^{+2}$ & $4$ &\\
    & $\Delta t_4$ & $2.50 \times 10^{+9}$ & $4.71 \times 10^{+2}$ & $4$ &\\
\hline
(B) & $\Delta t_1$ & $2.48 \times 10^{+6}$ & $2.11 \times 10^{+1}$ & $8$ &\\
    & $\Delta t_2$ & $2.49 \times 10^{+6}$ & $4.71 \times 10^{+1}$ & $8$ &\\
    & $\Delta t_3$ & $2.49 \times 10^{+6}$ & $6.66 \times 10^{+1}$ & $8$ &\\
    & $\Delta t_4$ & $2.49 \times 10^{+9}$ & $1.49 \times 10^{+2}$ & $8$ &\\
\hline
(C) & $\Delta t_1$ & $2.50 \times 10^{+8}$ & $6.66 \times 10^{+1}$ & $6$ &\\
    & $\Delta t_2$ & $2.50 \times 10^{+8}$ & $1.49\times 10^{+2}$ & $6$ &\\
    & $\Delta t_3$ & $2.50 \times 10^{+8}$ & $2.11 \times 10^{+2}$ & $6$ &\\
    & $\Delta t_4$ & $2.50 \times 10^{+8}$ & $4.71 \times 10^{+2}$ & $6$ &\\
\hline
(D) & $\Delta t_1$ & $5.00 \times 10^{+6}$ & $9.42 \times 10^{+0}$ & $6$ &\\
    & $\Delta t_2$ & $5.00 \times 10^{+6}$ & $2.11 \times 10^{+1}$ & $6$ &\\
    & $\Delta t_3$ & $5.00 \times 10^{+6}$ & $2.98 \times 10^{+1}$ & $6$ &\\
    & $\Delta t_4$ & $5.00 \times 10^{+6}$ & $6.66 \times 10^{+1}$ & $6$ &\\
\hline
\end{tabular}
\end{center}
}
\caption{Optimal parameters $\alpha_f$, $\alpha_p$ and number of stationary iterations to reach a tolerance of $10^{-8}$. On the left Table $\Delta t=0.01$, on the right Table $h=1/16$.}\label{table:testperiodic}
\end{table}

\begin{figure}
\includegraphics[scale=0.33]{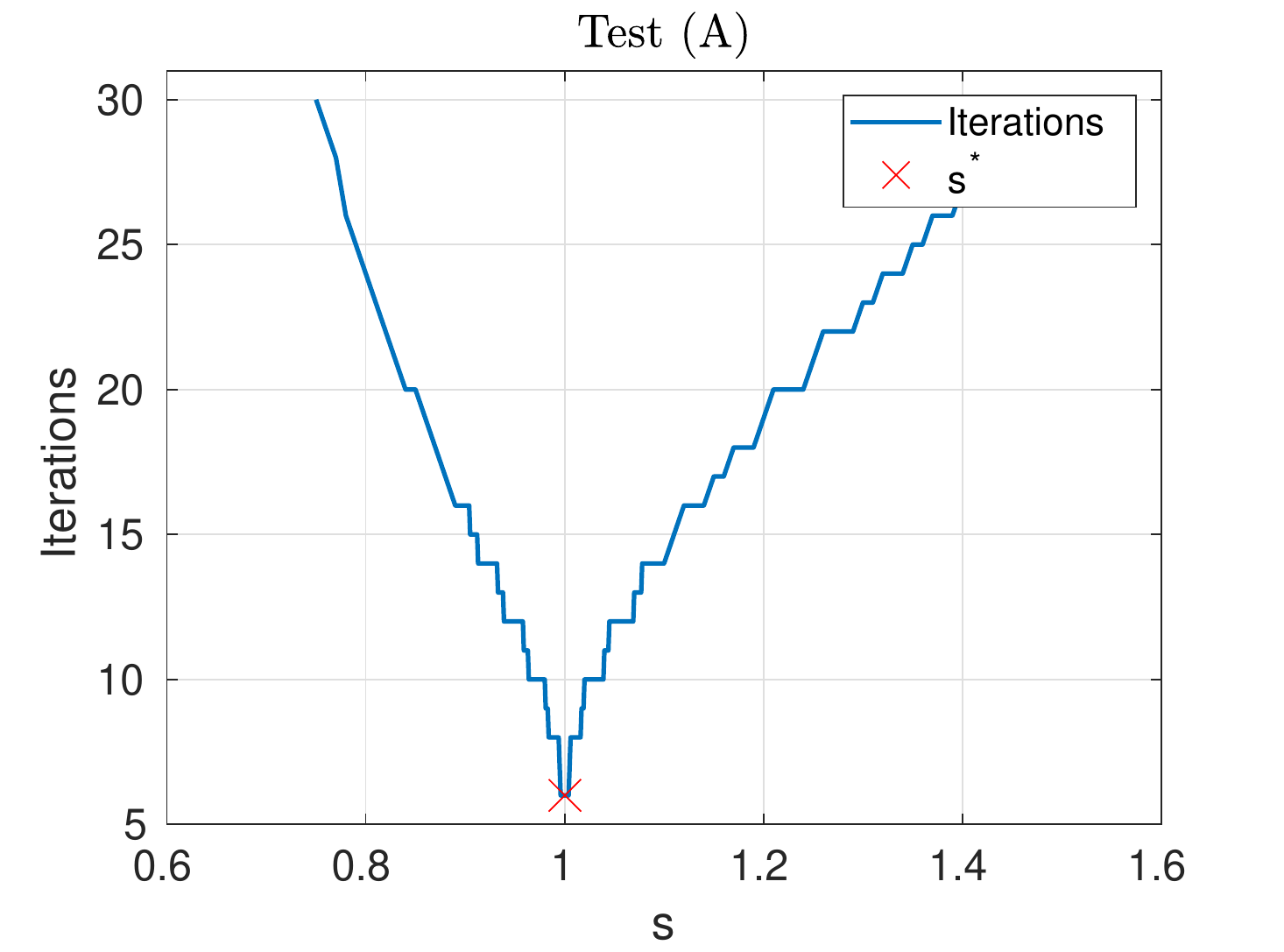}
\includegraphics[scale=0.33]{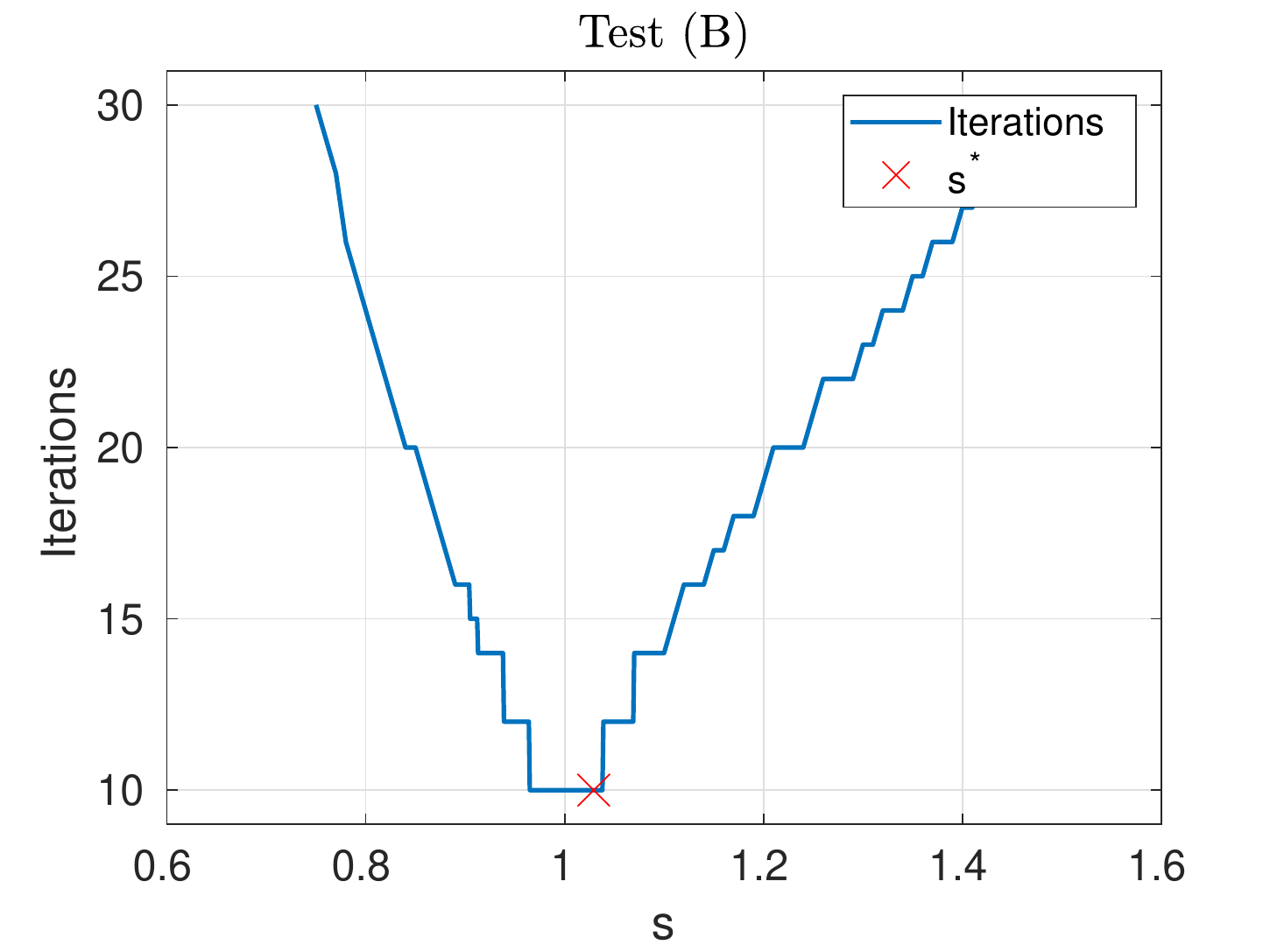}
\includegraphics[scale=0.33]{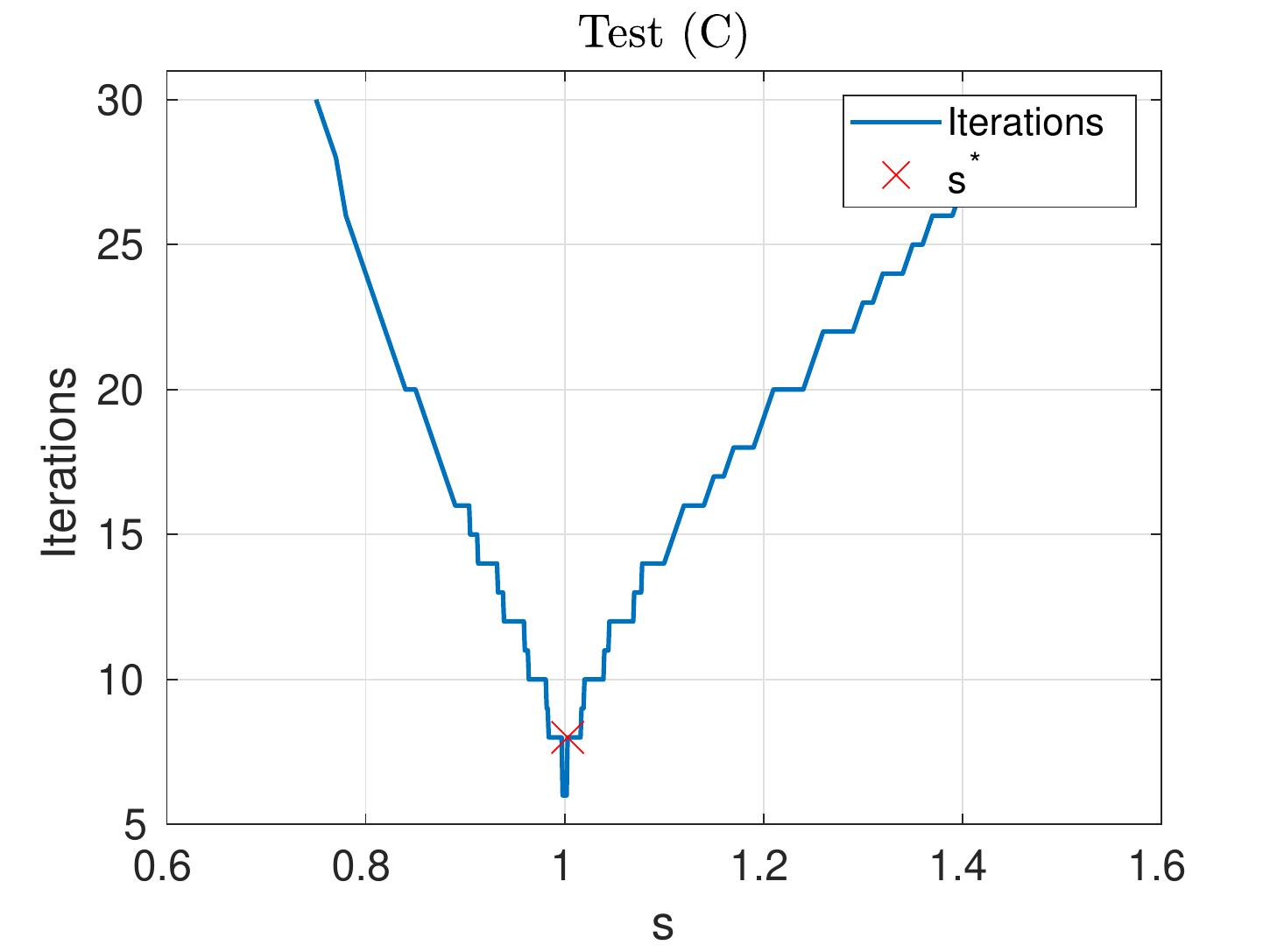}
\caption{Number of iterations to reach a relative tolerance of $10^{-8}$.}\label{fig:test_theovsnum}
\end{figure}

Figure \ref{fig:test_theovsnum} reports the number of iterations to solve the coupled system \eqref{eq:stationary_iteration} up to a relative tolerance of $10^{-8}$ for tests (A), (B) and (C). The discretization parameters are $h=1/32$, $\Delta=0.01$ and $\theta=0.5$. 
Notice that the optimal theoretical parameter $s^\star$ always leads to an optimal numerical convergence.

\subsection{Tests with analytic solution}\label{sec:test_analytic}

Consider the domain $\Omega_f = (0,0.5)\times(1,1.5)$, $\Omega_p = (0,0.5)\times(0.5,1)$ with interface $\Gamma = (0,1) \times \{1\}$, and the time interval $[0,0.5]$. The boundary conditions and forces $\mathbf{f}_f$ and $f_p$ are chosen in such a way that the exact solution of the problem is $\mathbf{u}_f = ( \sqrt{\mu_f \eta_p} \, \cos(t) , \, \alpha_{BJ} x \, \cos(t))$, $p_f = (2\mu_f (x+y-1) + (3\eta_p)^{-1})\, \cos(t)$, $p_p = (( -\alpha_{BJ} x(y-1)+ y^3/3 -y^2 +y )/\eta_p + 2 \mu_f x)\, \cos(t)$. In this test case, the exact solution satisfies the BJS condition \eqref{eq:interf3} at $\Gamma$ instead of the simplified condition \eqref{eq:interf3Simplified}.

For the space discretisation we consider $\mathbb{Q}_2-\mathbb{Q}_1$ elements for Stokes and $\mathbb{Q}_2$ elements for Darcy on uniform, structured computational grids made by rectangles, while for the time discretisation we use the backward Euler method ($\theta = 1$).

First, we test the robustness of the method with respect to the mesh size. For this purpose, we consider four computational meshes with sizes $h_j=0.1\times 2^{1-j}$, $j=1,\ldots,4$, and $\Delta t = 0.01$.
Then, we test the behaviour of the method with respect to $\Delta t$. To this aim, we consider a computational grid with $h=0.02$ and four timesteps $\Delta t_1 = 0.05$, $\Delta t_2 = 0.01$, $\Delta t_3 = 0.005$ and $\Delta t_4 = 0.001$ as in Sect. \ref{sec:periodic}.
In all cases, we solve the interface system \eqref{eq:interfaceSystem_v2} using GMRES \cite{Saad:1986:GMR} with tolerance $10^{-8}$ on the relative residual starting the iterations from ${\lambda}_f^{0} = 0$ and ${\lambda}_p^{0} = 0$ at the first time step, while at $t^n$ ($n \geq 1$) we set $\lambda_f^{n,0} = \lambda_f^{n-1}$ and $\lambda_p^{n,0} = \lambda_p^{n-1}$. The optimal value of $s$ in \eqref{eq:min_max} is computed using the same algorithm as in Sect. \ref{sec:periodic}.

Tables \ref{table:analyticH} and \ref{table:analyticDt} report the values of the optimal parameters $\alpha_f$ and $\alpha_p$ computed for the test cases (A)--(D) and the various discretization parameters, together with the number of iterations needed to solve \eqref{eq:interfaceSystem_v2} at $t^1$ and the rounded average number of iterations performed at successive time steps $t^n$ ($n>1$). The results show that the method is robust with respect to both the discretization parameters $h$ and $\Delta t$ and to the physical parameters that characterize the time-dependent Stokes-Darcy problem. As already observed in Sect. \ref{sec:periodic}, also in this case we notice that $\alpha_f$ is always few orders of magnitude larger than $\alpha_p$ so that the iterative method behaves like a Dirichlet-Robin one.

\begin{table}[bht]
\begin{center}
\begin{tabular}{ccllcc}
Test & Mesh & \multicolumn{1}{c}{$\alpha_f$} & \multicolumn{1}{c}{$\alpha_p$} & iter $t^1$ & iter $t^n$ \\
\hline
(A) & $h_1$ & $4.73 \times 10^{+7}$ & $1.05 \times 10^{+2}$ & $4$ & $2$ \\
    & $h_2$ & $2.38 \times 10^{+7}$ & $1.05 \times 10^{+2}$ & $4$ & $2$ \\
    & $h_3$ & $1.20 \times 10^{+7}$ & $1.05 \times 10^{+2}$ & $4$ & $3$ \\
    & $h_4$ & $5.99 \times 10^{+6}$ & $1.05 \times 10^{+2}$ & $4$ & $4$ \\
\hline
(B) & $h_1$ & $4.60 \times 10^{+4}$ & $3.33 \times 10^{+1}$ & $5$ & $4$ \\
    & $h_2$ & $2.35 \times 10^{+4}$ & $3.33 \times 10^{+1}$ & $6$ & $4$ \\
    & $h_3$ & $1.20 \times 10^{+4}$ & $3.33 \times 10^{+1}$ & $6$ & $4$ \\
    & $h_4$ & $6.04 \times 10^{+3}$ & $3.33 \times 10^{+1}$ & $8$ & $5$ \\
\hline
(C) & $h_1$ & $4.73 \times 10^{+6}$ & $1.05 \times 10^{+2}$ & $4$ & $3$ \\
    & $h_2$ & $2.38 \times 10^{+6}$ & $1.05 \times 10^{+2}$ & $4$ & $4$ \\
    & $h_3$ & $1.20 \times 10^{+6}$ & $1.05 \times 10^{+2}$ & $4$ & $4$ \\
    & $h_4$ & $6.00 \times 10^{+5}$ & $1.05 \times 10^{+2}$ & $6$ & $4$ \\
\hline
(D) & $h_1$ & $1.17 \times 10^{+5}$ & $1.49 \times 10^{+1}$ & $4$ & $4$ \\
    & $h_2$ & $4.69 \times 10^{+4}$ & $1.49 \times 10^{+1}$ & $4$ & $4$ \\
    & $h_3$ & $2.34 \times 10^{+4}$ & $1.49 \times 10^{+1}$ & $6$ & $4$ \\
    & $h_4$ & $1.19 \times 10^{+4}$ & $1.49 \times 10^{+1}$ & $6$ & $4$ \\
\hline
\end{tabular}
\caption{Optimal parameters $\alpha_f$, $\alpha_p$ and number of GMRES iterations for four computational meshes and fixed $\Delta t = 0.01$.}
\label{table:analyticH}
\end{center}
\end{table}

\begin{table}[bht]
\begin{center}
\begin{tabular}{ccllcc}
Test & $\Delta t$ & \multicolumn{1}{c}{$\alpha_f$} & \multicolumn{1}{c}{$\alpha_p$} & iter $t^1$ & iter $t^n$ \\
\hline
(A) & $\Delta t_1$ & $9.59 \times 10^{+6}$ & $4.71 \times 10^{+1}$ & $4$ & $4$ \\
    & $\Delta t_2$ & $9.57 \times 10^{+6}$ & $1.05 \times 10^{+2}$ & $4$ & $3$ \\
    & $\Delta t_3$ & $9.56 \times 10^{+6}$ & $1.49 \times 10^{+2}$ & $4$ & $3$ \\
    & $\Delta t_4$ & $9.51 \times 10^{+6}$ & $3.33 \times 10^{+2}$ & $4$ & $2$ \\
\hline
(B) & $\Delta t_1$ & $9.61 \times 10^{+3}$ & $1.49 \times 10^{+1}$ & $6$ & $5$ \\
    & $\Delta t_2$ & $9.56 \times 10^{+3}$ & $3.33 \times 10^{+1}$ & $6$ & $4$ \\
    & $\Delta t_3$ & $9.52 \times 10^{+3}$ & $4.71 \times 10^{+1}$ & $7$ & $4$ \\
    & $\Delta t_4$ & $9.40 \times 10^{+3}$ & $1.05 \times 10^{+2}$ & $8$ & $4$ \\
\hline
(C) & $\Delta t_1$ & $9.60 \times 10^{+5}$ & $4.71 \times 10^{+1}$ & $5$ & $4$ \\
    & $\Delta t_2$ & $9.58 \times 10^{+5}$ & $1.05 \times 10^{+2}$ & $5$ & $4$ \\
    & $\Delta t_3$ & $9.57 \times 10^{+5}$ & $1.49 \times 10^{+2}$ & $5$ & $4$ \\
    & $\Delta t_4$ & $9.51 \times 10^{+5}$ & $3.33 \times 10^{+2}$ & $6$ & $3$ \\
\hline
(D) & $\Delta t_1$ & $1.90 \times 10^{+4}$ & $6.66 \times 10^{+0}$ & $6$ & $4$ \\
    & $\Delta t_2$ & $1.88 \times 10^{+4}$ & $1.49 \times 10^{+1}$ & $6$ & $4$ \\
    & $\Delta t_3$ & $1.87 \times 10^{+4}$ & $2.11 \times 10^{+1}$ & $6$ & $4$ \\
    & $\Delta t_4$ & $3.91 \times 10^{+4}$ & $4.71 \times 10^{+1}$ & $6$ & $4$ \\
\hline
\end{tabular}
\caption{Optimal parameters $\alpha_f$, $\alpha_p$ and number of GMRES iterations for four values of $\Delta t$ and fixed $h = 0.02$.}
\label{table:analyticDt}
\end{center}
\end{table}

%
%

\subsection{Tests without analytic solution}\label{sec:realistic}

Consider the dimensionless domain $\Omega = (0,1)\times (0,1)$ with interface at $y=0.4$ and $\Omega_f$ the upper subdomain, and the dimensionless time interval $[0,t_f]$ with $t_f=1$. For Darcy's problem, we impose homogeneous Dirichlet boundary condition at the bottom boundary $(0,1)\times\{0\}$ of the domain, and homogeneous Neumann boundary conditions on the remaining external boundaries of $\Omega_p$. In $\Omega_f$, we impose homogeneous Dirichlet boundary conditions on the external lateral boundaries, while we set $\mathbf{u}_f = (u_1,0)$ on the top boundary $(0,1)\times\{1\}$ with $u_1 = \min (2U_f\,t/t_f,U_f)$ and dimensionless velocity $U_f=1$. The forces $\mathbf{f}_f$ and $f_p$ are both zero. The physical parameters are chosen as in tests (A) and (D) of Sect. \ref{sec:choice_parameters}, see Table \ref{table:values}. 

\begin{figure}[bht]
\begin{center}
\includegraphics[width=0.45\textwidth]{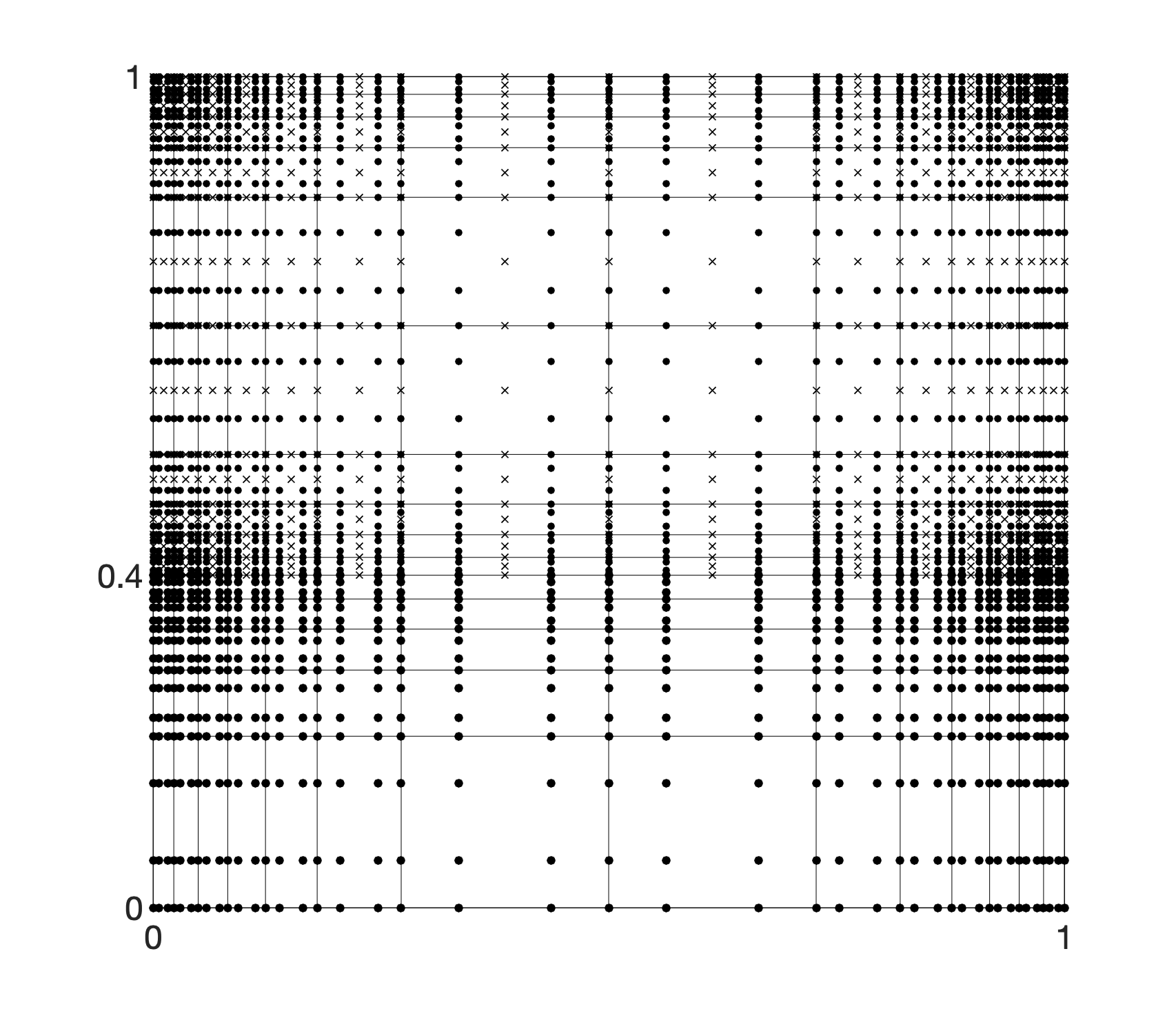}
\includegraphics[width=0.45\textwidth]{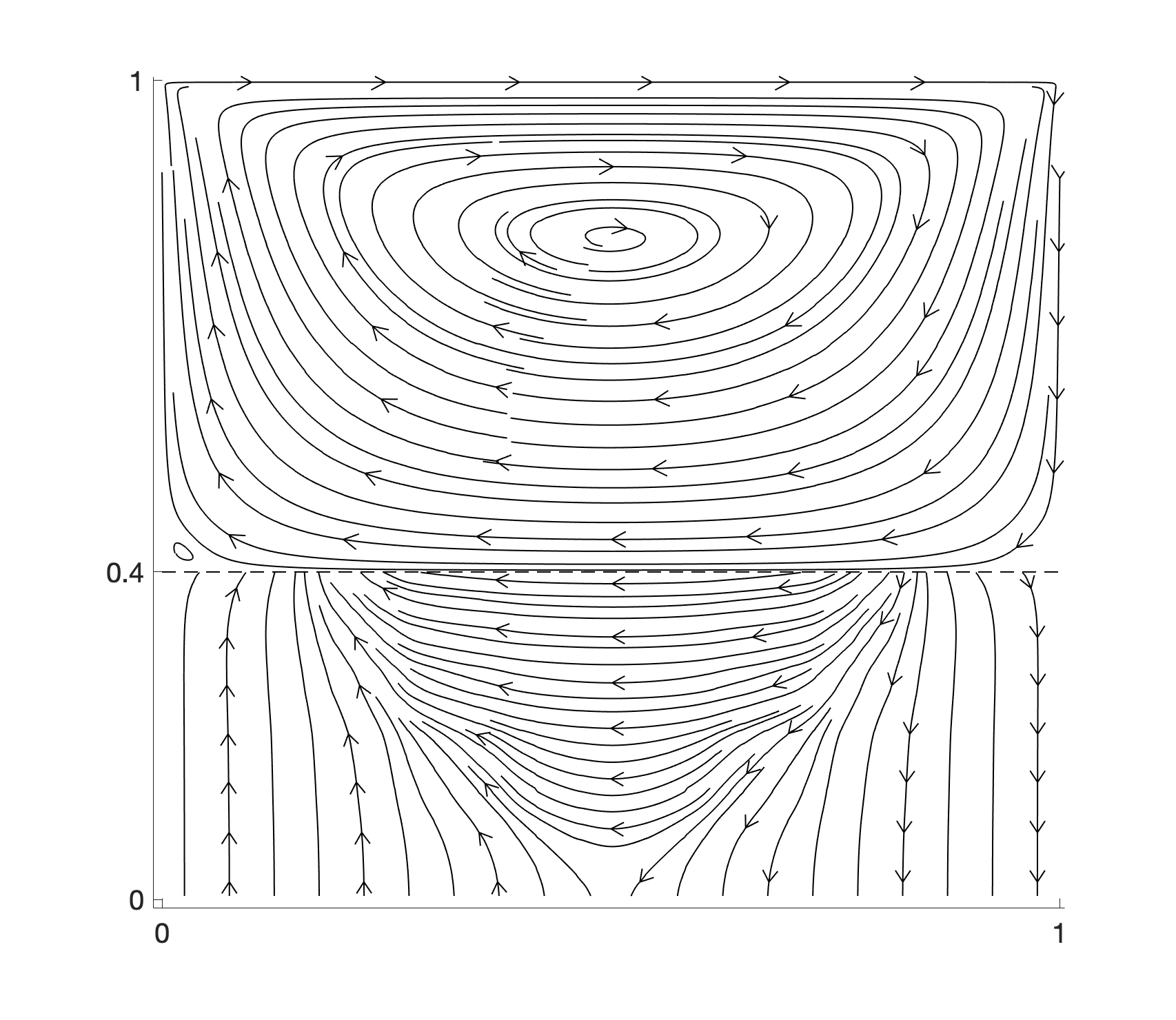}
\caption{Sample computational mesh corresponding to $h=h_1$ (left) and representation of the streamlines of the computed solution at $t=0.75$ using mesh with $h=h_2$.}\label{fig:meshSolNoAnalytic}.
\end{center}
\end{figure}

For the discretization, we use $\mathbb{Q}_3-\mathbb{Q}_2$ elements for Stokes and $\mathbb{Q}_3$ elements for Darcy with Gauss-Lobatto nodes for the $\mathbb{Q}_3$ polynomials, and the Crank-Nicolson method ($\theta=\frac{1}{2}$). The meshes are non-uniform with smaller elements in the neighbourhood of the external boundary and of the interface. An example is shown in Fig. \ref{fig:meshSolNoAnalytic} (left).

Tables \ref{table:noAnalyticH} and \ref{table:noAnalyticDt} indicate the values of the optimal parameters $\alpha_f$ and $\alpha_p$ and the number of iterations needed to solve \eqref{eq:interfaceSystem_v2} at $t^1$ and the rounded average number of iterations performed at successive time steps $t^n$ ($n>1$). For the computation of the optimal parameters, an average value of $h$ at the interface is used due to the non-uniformity of the mesh. In Table \ref{table:noAnalyticH}, we set $\Delta t = 0.05$ and we consider three values of $h$: $h_1 \approx 0.0714$, $h_2 \approx 0.0357$, $h_3 = 0.0179$. In Table \ref{table:noAnalyticDt}, we consider the mesh characterized by $h_2$ and we consider three values of $\Delta t$: $\Delta t_1 = 0.125$, $\Delta t_2 = 0.0625$, $\Delta t_3 = 0.03125$. As already observed in Sect. \ref{sec:test_analytic}, the method is robust with respect to both the discretization and the physical parameters and the choice of the optimal parameters $\alpha_f$ and $\alpha_p$ indicate a behaviour analogous to a Dirichlet-Robin algorithm.

\begin{figure}[bht]
\includegraphics[width=0.45\textwidth]{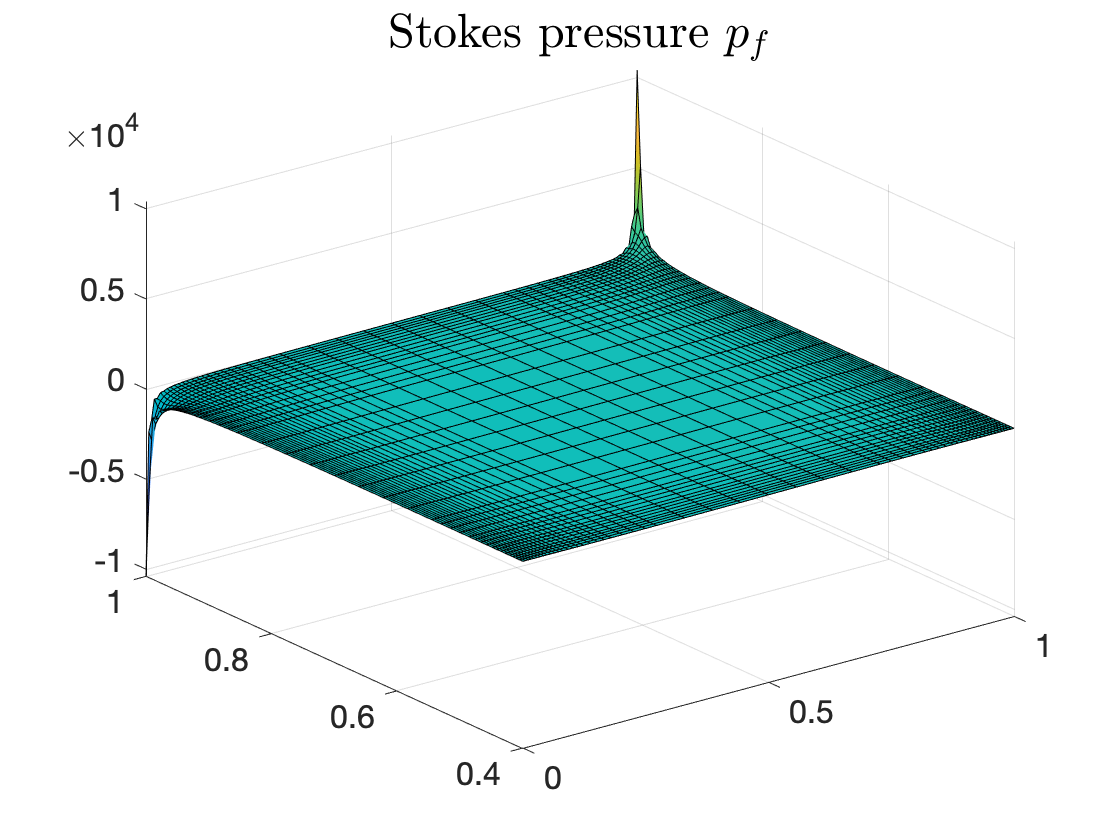}
\includegraphics[width=0.45\textwidth]{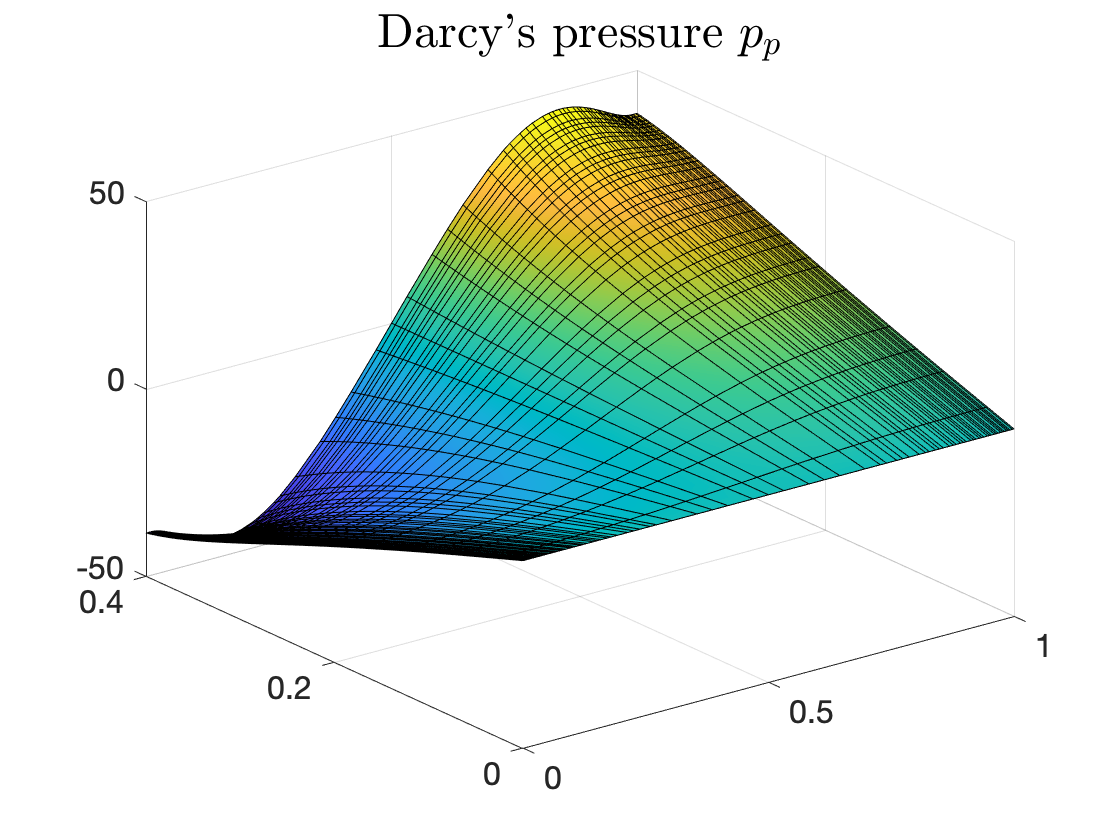}
\caption{Stokes pressure $p_f$ (left) and Darcy's pressure $p_p$ (right) computed at $t=0.75$ using mesh with $h=h_2$.}\label{fig:pressureNoAnalytic}.
\end{figure}

Figure \ref{fig:meshSolNoAnalytic} (right) and Fig. \ref{fig:pressureNoAnalytic} show the velocity and pressure computed at $t=0.75$ using the mesh characterized by size $h_2$ and $\Delta t = 0.05$. (The velocity in Darcy's domain has been postprocessed from the pressure $p_p$ using the MATLAB command `gradient'.)

\begin{table}[bht]
\begin{center}
\begin{tabular}{cccllcc}
Test & Mesh & No. unknowns at $\Gamma$ & \multicolumn{1}{c}{$\alpha_f$} & \multicolumn{1}{c}{$\alpha_p$} & iter $t^1$ & iter $t^n$ \\
\hline
(A) & $h_1$ &  86 & $2.28 \times 10^{+7}$ & $3.33 \times 10^{+1}$ & $12$ & $9$ \\
    & $h_2$ & 170 & $1.14 \times 10^{+7}$ & $3.33 \times 10^{+1}$ & $12$ & $9$ \\
    & $h_3$ & 338 & $5.71 \times 10^{+6}$ & $3.33 \times 10^{+1}$ & $14$ & $10$ \\
\hline
(D) & $h_1$ &  86 & $4.43 \times 10^{+4}$ & $4.71$ & $8$ & $7$ \\
    & $h_2$ & 170 & $2.25 \times 10^{+4}$ & $4.71$ & $10$ & $8$ \\
    & $h_3$ & 338 & $1.13 \times 10^{+4}$ & $4.71$ & $12$ & $9$ \\
\hline
\end{tabular}
\caption{Test of Sect. \ref{sec:realistic}: optimal parameters $\alpha_f$, $\alpha_p$ and number of GMRES iterations for three non-uniform computational meshes and $\Delta t = 0.05$.}
\label{table:noAnalyticH}
\end{center}
\end{table}

\begin{table}[bht]
\begin{center}
\begin{tabular}{ccllcc}
Test & $\Delta t$ & \multicolumn{1}{c}{$\alpha_f$} & \multicolumn{1}{c}{$\alpha_p$} & iter $t^1$ & iter $t^n$ \\
\hline
(A) & $\Delta t_1$ & $1.14 \times 10^{+7}$ & $2.11 \times 10^{+1}$ & $12$ & $9$ \\
    & $\Delta t_2$ & $1.14 \times 10^{+7}$ & $2.98 \times 10^{+1}$ & $12$ & $9$ \\
    & $\Delta t_3$ & $1.14 \times 10^{+7}$ & $4.21 \times 10^{+1}$ & $12$ & $9$ \\
\hline  
(D) & $\Delta t_1$ & $2.26 \times 10^{+4}$ & $2.98$ & $10$ & $8$ \\
    & $\Delta t_2$ & $2.25 \times 10^{+4}$ & $4.21$ & $10$ & $8$ \\
    & $\Delta t_3$ & $2.24 \times 10^{+4}$ & $5.96$ & $10$ & $8$ \\
\hline
\end{tabular}
\caption{Test of Sect. \ref{sec:realistic}: optimal parameters $\alpha_f$, $\alpha_p$ and number of GMRES iterations for three values of $\Delta t$ and mesh with $h_2$.}
\label{table:noAnalyticDt}
\end{center}
\end{table}

\section{Conclusions}

In this paper, we formulated and analyzed an optimized Schwarz method for the time-dependent Stokes-Darcy problem. Since the convergence factor is different from other cases studied in the literature, we proposed a novel approach to characterize the optimal parameters in the interface Robin conditions in order to guarantee robustness of the iterative method with respect to physical and discretization parameters. Numerical experiments carried out for various configurations of the problem showed the effectiveness of the studied method.

\bigskip

\paragraph{Acknowledgements.} The first author acknowledges funding by the EPSRC grant EP/V027603/1.

\bibliographystyle{plain}
\bibliography{DV}

\end{document}